\documentclass[12pt]{article}
\usepackage{amsfonts,amssymb,amsthm,fullpage}
\title{The congruence kernel of an arithmetic lattice in a rank
one algebraic group over a local field}
\author{A. W. Mason$^1$, A. Premet$^{2}$, B. Sury$^{3}$ \footnote{Partially supported
by an EPSRC Visiting Fellowship GR/N32211/01.} , P. A.
Zalesskii$^{4}$ \footnote{Partially supported by the Edinburgh
Mathematical Society Research Support Fund, the Glasgow
Mathematical Journal Trust Fund, Conselho Nacional de
Desenvolvimento Cient\'{i}fico e Tecnol\'{o}gico (CNPq) and
Coordena\c{c}\~{a}o de Aperfei\c{c}oamento de Pessoal de N\'{i}vel
Superior (Capes).}}
\date{}

\def\char{\mathop{\mathrm{char}}}

\def\mod{\mathop{\mathrm{mod}}}
\newtheorem{thm}{Theorem}[section]
\newtheorem{lem}[thm]{Lemma}
\newtheorem{cor}[thm]{Corollary}
\begin{document}
 \maketitle
\noindent $^1$ Department of Mathematics, University of Glasgow,
Glasgow G12 8QW, Scotland UK (e-mail: awm@maths.gla.ac.uk)\\
\noindent $^2$ School of Mathematics, University of Manchester,
Oxford Road, Manchester M13 9PL, UK (e-mail:
sashap@maths.man.ac.uk)\\
\noindent $^3$ Statistics-Mathematics Unit, Indian Statistical
Institute, Bangalore 560 059, India (e-mail: sury@isibang.ac.in)\\
\noindent $^4$ Department of Mathematics, University of Brasilia,
70.910 Brasilia DF, Brazil (e-mail: pz@mat.unb.br)
\section*{Abstract}
Let $k$ be a global field and let $k_v$ be the completion of $k$
with respect to $v$, a non-archimedean place of $k$. Let
$\mathbf{G}$ be a connected, simply-connected algebraic group over
$k$, which is absolutely almost simple of $k_v$-rank 1. Let
$G=\mathbf{G}(k_v)$. Let $\Gamma$ be an {\it arithmetic lattice}
in $G$ and let $C=C(\Gamma)$ be its {\it congruence kernel}.
Lubotzky has shown that $C$ is infinite, confirming an earlier
conjecture of Serre. Here we provide complete solution of the {\it
congruence subgroup problem} for $\Gamma$ by determining the {\it
structure} of $C$. It is shown that $C$ is a free profinite
product, one of whose factors is $\hat{F}_{\omega}$, the free
profinite group on countably many generators. The most surprising
conclusion from our results is that the structure of $C$ depends
{\it only} on the characteristic of $k$. The structure of $C$ is
already known for a number of special cases. Perhaps the most
important of these is the ({\it non-uniform}) example $\Gamma={\rm
SL}_2(\mathcal{O}(S))$, where $\mathcal{O}(S)$ is the {\it the
ring of S-integers in k}, with $S=\{v\}$, which plays a central
role in the theory of Drinfeld modules. The proof makes use of a
decomposition theorem of Lubotzky, arising from the action of
$\Gamma$ on the {\it Bruhat-Tits tree} associated with $G$.\\

\noindent {\it 2000 Mathematics Subject Classification} : 20G30,
11F06, 20E08, 20E18.\\

\section*{Introduction}

Let $k$ be a global field and let $\mathbf{G}$ be a connected,
simply-connected linear algebraic group over $k$, which is
absolutely almost simple. For each non-empty, finite set $S$ of
places of $k$, containing all the archimedean places, let
$\mathcal{O}(S)$ denote the corresponding {\it ring of} $S$-{\it
integers in} $k$. The problem of determining whether or not a
finite index subgroup of the arithmetic group,
$\mathbf{G}(\mathcal{O}(S))$, contains a principal congruence
subgroup (modulo some non-zero $\mathcal{O}(S)$-ideal), the
so-called {\it congruence subgroup problem} or CSP, has attracted
a great deal of attention since the $19$th century. As a measure
of the extent of those finite index subgroups of
$\mathbf{G}(\mathcal{O}(S))$ which are not congruence, its
so-called {\it non-congruence subgroups}, Serre [S1] has
introduced a profinite group, $C(S,\mathbf{G})$, called the {\it
(S-)congruence kernel} of $\mathbf{G}$. In his terminology [S1]
the CSP for this group has an {\it affirmative} answer if this
kernel is finite. Otherwise the CSP has an {\it essentially
negative} answer. The principal result in [S1] is that, for the
case $\mathbf{G}=\mathbf{SL}_2$, the congruence kernel
$C(S,\mathbf{G})$ is {\it finite} if and only if card$S\geq2$.
Moreover Serre has formulated the famous {\it congruence subgroup
conjecture} [PR, p.~556], which states that the answer to the CSP
is determined entirely by the $S$-$rank$ of $\bf{G}$,  $
\rm{rank}_S\bf{G}$. (See [Mar, p.~258].) It is known [Mar, (2.16)
Theorem, p.~269] that $C(S,\bf{G})$ is finite (cyclic), when
$\bf{G}$ is $k$-$isotropic$ and rank$_S\bf{G}\geq2$. It is also
known that $C(S,\bf{G})$ is infinite for many ``rank one" $\bf{G}$
(for example, $\bf{G}=\bf{SL}_2$). The conjecture however remains
open for some of these cases. (See, for example, [L3].) The
congruence kernel $C(S,H)$ can be defined in a similar way for
every subgroup $H$ of $\mathbf{G}(k)$ which is commensurable with
$\bf{G}(\cal{O}(S))$. (From this definition it is
clear that $C(S,H)$ is finite if and only if $C(S,\bf{G})$ is finite.)\\

The books of Margulis [Mar, p.~268] and Platonov/Rapinchuk [PR,
Section 9.5] emphasise the importance of determining the {\it
structure} of the congruence kernel. (Lubotzky refers to this as
the {\it complete} solution of the CSP.) In this paper we are
concerned with the structure of infinite congruence kernels. The
first result of this type is due to Mel'nikov [Me], who shows
that, for the case where $\mathbf{G}=\mathbf{SL}_2$,
$k=\mathbb{Q}$ and $S=\{\infty\}$, (i.e.
$\mathbf{G}(\mathcal{O}(S)) = { \rm SL}_2(\mathbb{Z})$, the
classical {\it modular group}), the congruence kernel is
isomorphic to $\hat{F}_{\omega}$, the {\it free profinite group on
countably many generators}. Lubotzky [L1] has proved that, when
$\bf{G}=\bf{SL}_2$ and card $S=1$, the congruence kernel of ${ \rm
SL}_2(\cal{O}(S))$ has a closed subgroup isomorphic to
$\hat{F}_{\omega}$, reproving Mel'nikov's result in the process.
(When char $k=0$ and card $S=1$,
 it is known that $ k= \mathbb{Q}$ or $\mathbb{Q}(\sqrt{-d})$, with $S=\{ \infty\}$, where $d$
is a square-free positive rational integer.) In [Mas2] it is shown
that, when $\mathbf{G}=\mathbf{SL}_2$ and card $S=1$,
 the congruence kernel maps onto every free profinite group of finite rank.\\

In this paper we use different methods to determine the structure of
the congruence kernel of an arithmetic lattice in a rank one
algebraic group over a local field, providing a complete solution of
the CSP for this case. With the above notation let $V_k$ be the set
of places of $k$ and let (the local field) $k_v$
 be the completion of $k$ with respect to $v$. In addition to the
above hypotheses we assume that $\mathbf{G}$ has $k_v$-rank 1. We
denote the set of $k_v$-rational points, $\mathbf{G}(k_v)$, by
$G$. Let $\Gamma$ be a {\it lattice} in $G$, i.e. a discrete
subgroup of (the locally compact group) $G$ for which
$\mu(G/\Gamma)$ is {\it finite}, where $\mu$ is a Haar measure on
$G$. As usual $\Gamma$ is said to be {\it cocompact} (resp. {\it
non-uniform}) if $G/\Gamma$ is compact (resp. not compact). We
assume further that $\Gamma$ is $(S-)${\it arithmetic}, i.e.
$\Gamma$ is commensurable with $\mathbf{G}(\mathcal{O})$, where
$\cal{O}=\cal{O}(S)$ is as above.\\

\noindent {\it Example}.  When $\char k > 0$, $S = \{v\}$ and
$\mathbf{G}=\mathbf{SL}_2$, the group $\Gamma = {\rm
SL}_2(\cal{O})$ is a (non-uniform) arithmetic lattice (in ${\rm
SL}_2(k_v)$). This lattice, which plays a central role in the
theory of Drinfeld modules, is
the principal focus of attention in Chapter II of Serre's book [S2]. \\

As in Margulis's book [Mar, Chapter I, 3.1, p.60] we assume that
$\mathbf{G}$ is $k$-subgroup of $\mathbf{GL}_n$, for some $n$. We
consider the standard representation for $\mathbf{GL}_n(k_v)$. For
each $\cal{O}$-ideal $\mathfrak{q}$, we put
$$
\mathbf{GL}_n(\mathfrak{q}) = \left\{X \in
\mathbf{GL}_n(\mathcal{O})\,|\,\, X \equiv I_n\:(\mod
\frak{q})\right\}.
$$

\noindent We denote $\mathbf{G}\cap\mathbf{GL}_n(\mathfrak{q})$, the
{\it principal S-congruence subgroup of $\mathbf{G}$ (of level
$\mathfrak{q}$)}, by $\mathbf{G}(\mathfrak{q})$. If $M$ is any
subgroup of $G$ commensurable with $\mathbf{G}(\mathcal{O})$ we put
$M(\mathfrak{q})=M\cap\mathbf{G}(\mathfrak{q})$. It is clear that
$M(\mathfrak{q})$ is of finite index in $M$ when
$\mathfrak{q}\neq\{0\}$. (We note that although the definition of
$\mathbf{G}(\mathcal{O}(S))$ depends on the $k$-embedding of
$\mathbf{G}$ into $\mathbf{GL}_n$, the class of the $S$-arithmetic
subgroups does not.)

\noindent The finite index subgroups of $\Gamma(\mathcal{O})$
define the {\it S-arithmetic} topology on $\Gamma$. The completion
of $\Gamma$ with respect to this topology is a profinite group
denoted by $\hat{\Gamma}$. On the other hand the subgroups
$\Gamma(\mathfrak{q})$, where $\mathfrak{q}\neq\{0\}$, define the
{\it S-congruence} topology on $\Gamma$ and the completion of
$\Gamma$ with respect to this topology is also a profinite group
denoted by $\bar{\Gamma}$. Since every $S$-congruence subgroup is
$S$-arithmetic, there is an exact sequence
$$
1 \rightarrow C(\Gamma) \rightarrow \hat{\Gamma} \rightarrow
\overline{\Gamma} \rightarrow 1.
$$

\noindent The (profinite) group $C(\Gamma)(=C(S,\Gamma))$ is
called the {\it (S-)congruence kernel} of $\Gamma$. It is known
[Mar Chapter I, 3.1] that the definition of $C(\Gamma)$ does not
depend on the choice of $k$-representation. (The definition of
congruence kernel extends to any $S$-arithmetic subgroup of $G$,
including any finite index subgroup of $\Gamma$.)\\ \noindent Our
principal results are the following.

\medskip

\noindent{\bf Theorem A.}  {\it If  $\Gamma$ is cocompact, then}
$$
C(\Gamma)\cong\hat{F}_{\omega}.
$$

\smallskip

It is well-known that $\Gamma$ is cocompact when, for example, char
$k=0$. Here Theorem A applies to the case where $S$ consists of
precisely one non-archimedean place, together with all the
archimedean places, and $\mathbf{G}$ is anisotropic over all the
archimedean places. For examples of cocompact lattices of the above
type in ${\rm SL}_2(\mathbb{Q}_p)$, where $\mathbb{Q}_p$ is the
$p$-adic completion of $\mathbb{Q}$,  see [S2, p.~84]. \noindent
This result however is not a straightforward generalization of
Mel'nikov's theorem [Me]. On the one hand ${ \rm SL}_2(\mathbb{Z})$
is {\it not} a lattice in ${\rm SL}_2(\mathbb{Q}_p)$. On the other
hand ${\rm SL}_2(\mathbb{Z})$ is a {\it non-uniform} lattice in
${\rm SL}_2(\mathbb{R})$. (See [Mar, p.~295].) Moreover the fourth
author [Za2] has proved that the congruence kernel of every
arithmetic lattice in ${\rm SL}_2(\mathbb{R})$ is isomorphic to
$\hat{F}_{\omega}$. Lattices to which Theorem A refers have a free,
non-cyclic subgroup of finite index. (See Lemma 2.1) Consequently
this result does {\it not} apply to the Bianchi groups, ${\rm
SL}_2(\mathcal{O}_d)$, where $\mathcal{O}_d$ is the ring of integers
in the imaginary quadratic number field $\mathbb{Q}(\sqrt{-d})$,
with $d >0$.
\medskip

\noindent{\bf Theorem B.} {\it If  $\Gamma$ is non-uniform and
$p=$ {\rm char} $k$, then
$$
C(\Gamma)\cong \hat{F}_{\omega}\amalg N(\Gamma),
$$

\smallskip

\noindent
the free profinite product of $\hat{F}_{\omega}$ and
$N(\Gamma)$, where $N(\Gamma)$ is a free profinite product of
groups, each of which is isomorphic to the direct product of
$2^{{\aleph}_0}$
copies of $\mathbb{Z}/p\mathbb{Z}$.}\\

\noindent The most interesting consequence of Theorems A and B is
that the structure of $C(\Gamma)$ depends {\it only} on the
characteristic of $k$.\\ \noindent The proofs are based on the
action of $G$, and hence $\Gamma$, on the associated Bruhat-Tits
tree $T$. The theory of groups acting on trees shows how to derive
the structure of $\Gamma$ from that of the quotient graph
$\Gamma\backslash T$. For the cocompact case it is well known that
$\Gamma\backslash T$ is finite. Theorem A then follows
from the theory of free profinite groups.\\
\noindent For the non-uniform case the situation is much more
complicated. Here Lubotzky [L2] has shown that $\Gamma\backslash T$
is the union of a finite graph together with a (finite) number of
ends, each of which corresponds to $\mathbf{P}$, a minimal parabolic
$k_v$-subgroup of $\mathbf{G}$. The proof that the torsion-free part
of the decomposition of $C(\Gamma)$ is $\hat{F}_{\omega}$ involves
substantially more effort than that of Theorem A. It is shown that
the torsion part $N(\Gamma)$ is a free profinite product of groups
each isomorphic to $C(U)=C(\mathbf{U}(\mathcal{O}))$, the {\it
S-congruence kernel} of $\mathbf{U}$, where $\mathbf{U}$ is the
unipotent radical of some $\mathbf{P}$ of the above type. Unlike the
characteristic zero unipotent groups, which have trivial congruence
kernel, the congruence kernel $C(U)$ is huge in positive
characteristic. The various ends of the quotient graph correspond to
certain unipotent subgroups and their congruence kernels contribute
to $N(\Gamma)$. It is known [BT2] that such a $\mathbf{U}$, and
hence $C(U)$, is nilpotent of class at most $2$. In fact we show
that $C(U)$ is {\it abelian}, even when $\mathbf{U}$ is not. In the
proofs the various types of $\mathbf{G}$, which arise from Tits
Classification [T], are dealt with separately. A crucial ingredient
(when dealing with
non-abelian $\mathbf{U}$) is the following unexpected property of "rank one" unipotent radicals. \\

\medskip

\noindent {\bf Theorem C.} {\it Let $\mathbf{U}$ be the unipotent
radical of a minimal parabolic $k_v$-subgroup of $\mathbf{G}$ of
the above type (so that $\mathbf{U}$ is defined over $k_v$). If
$\mathbf{U}(k)$ is not
abelian then $\mathbf{U}$ is defined over $k$.}\\

\noindent For our purposes the importance of Theorem C is that it
ensures that the structure of $C(U)$ needs to be determined only for
one particular $\mathbf{U}$. Theorem B extends a number of existing
results. The fourth author [Za1, Theorem 4.3] has proved Theorem B
for the special case $\mathbf{G} =\mathbf{SL}_2$ and $S=\{v\}$.
(This case is rather more straightforward since here $\mathbf{U}$ is
abelian, and so Theorem C, for example, is not required.) Lubotzky
[L1] has proved that, for this case, $C(\Gamma)$ has a closed
subgroup isomorphic to $\hat{F}_{\omega}$. Lubotzky has also shown
[L2,
Theorem 7.5] that $C(\Gamma)$ is {\it infinite} when $\Gamma$ is non-uniform.\\

\noindent Let $\mathbf{H}$ be any semisimple algebraic group over
$k$. In addition to the $S$-congruence kernel, $C(S,\bf{H})$,
there is another group called the {\it S-metaplectic kernel},
$M(S, \bf{H})$, whose definition (originally due to Moore) is
cohomological. (See, for example, [PR, p. 557].) It is known [PRr,
Theorem 9.15, p.557] that these groups are closely related when
$C(S,\bf{H})$ is finite. (The structure of $M(S, \bf{H})$ has been
determined for many such cases; see [PRap].) In this paper however
we are concerned with {\it infinite} congruence kernels.

\section{ Arithmetic lattices}

This section is devoted to a number of properties of arithmetic
lattices which are needed to establish our principal results. From
now we will use {\it lattice} as an abbreviation for {\it lattice
in} $G=\mathbf{G}(k_v)$, where $\mathbf{G}$ and $k_v$ are defined as
above. We begin with a general property of lattices.

\begin{lem} If  $\Gamma$ is any lattice, then
 $\Gamma$ is not virtually solvable.
\end{lem}
 \noindent{\bf Proof.} It is known that $\Gamma$ is {\it Zariski-dense} in
 $\mathbf{G}$. (See [Mar, (4.4) Corollary, p.~93] and [Mar, (2.3)
 Lemma, p.~84].)
 \noindent It follows that $[\Gamma,\Gamma]$ is Zariski-dense in
 $[\mathbf{G},\mathbf{G}]=\mathbf{G}$, by [B, Proposition, p.~59]
 and [B, Proposition, p.~181]. If $\Gamma$ is virtually solvable
 then $\mathbf{G}$ is finite, which contradicts the fact that it
 has $k_v$-rank 1.$\hfill\Box$

 \medskip

For each non-archimedean $v \in V_k$, we denote the completion of
$\mathcal{O}$ with respect to $v$ by $\mathcal{O}_v$. This is a
local ring with a finite residue field. Recall that the {\it
restricted topological product} is defined as
$$
\mathbf{G}(\hat{\mathcal{O}}) = \prod_{v \not\in S}
\mathbf{G}(\mathcal{O}_v);
$$
see [PR, p.~161]. The group $\mathbf{G}(\hat{\mathcal{O}})$ is a
topological group with a base of neighbourhoods of the identity
consisting of all subgroups of the form
$$
\prod_{v \not\in S} M_v, \eqno{(*)}
$$
where each $M_v$ is an open subgroup of
$\mathbf{G}(\mathcal{O}_v)$ and $M_v = \mathbf{G}(\mathcal{O}_v)$,
for all but finitely many $v \not\in S$. Let $\mathfrak{m}$ denote
the maximal ideal of the (local) ring $\mathcal{O}_v$. Then the
"principal congruence subgroups", $\mathbf{G}(\mathfrak{m}^t)$,
where $t\geq 1$, provide a base of neighbourhoods of the identity
in $\mathbf{G}(\mathcal{O}_v)$; see [PR, p.~134]. The group
$\mathbf{G}(\mathcal{O})$ embeds, via the ``diagonal map'', in
$\mathbf{G}(\hat{\mathcal{O}})$. Let
$\overline{\mathbf{G}}(\mathcal{O})$ denote the "congruence
completion" of $\mathbf{G}(\mathcal{O})$ determined by its
$S$-congruence subgroups. The hypotheses on $\mathbf{G}$ ensure
that the following holds.

\begin{lem}``The Strong Approximation Property''
$$
\overline{\mathbf{G}}(\mathcal{O}) \cong
\mathbf{G}(\hat{\mathcal{O}}).
$$
\end{lem}

\noindent{\bf Proof.} By [PR, Theorem 7.12, p.~427] it suffices to
verify that
$$
G_S := \prod_{v\in S}\mathbf{G}(\mathcal{O}_v)
$$
is not compact. \noindent Now by [Mar, (3.2.5), p.~63] the group
$\mathbf{G}(\mathcal{O})$ is a lattice in $G_S$. If $G_S$ is
compact then $\mathbf{G}(\mathcal{O})$ and hence $\Gamma$ are
finite, which contradicts Lemma 1.1.$\hfill\Box$

\medskip

We record another well-known property of $\Gamma$.

\begin{lem}  With the above notation,
$$
C(\Gamma) = \displaystyle{\bigcap_{\frak{q} \neq \{0\}}}
\hat{\Gamma}(\frak{q}).
$$
\end{lem}

\noindent It follows that, for all $\mathfrak{q}\neq\{0\}$, there
is an exact sequence
$$
1 \rightarrow C(\Gamma) \rightarrow \hat{\Gamma}(\mathfrak{q})
\rightarrow \overline{\Gamma}(\mathfrak{q}) \rightarrow 1.
$$
\noindent More generally let $M$ be any group of matrices over
$\mathcal{O}$. For each non-zero $\mathcal{O}$-ideal
$\mathfrak{q}$ we define the (finite index) subgroup
$M(\mathfrak{q})$ of $M$ in the natural way as above. Then the
subgroups $M(\mathfrak{q})$ form a base of neighbourhoods of the
identity in $M$ for the {\it congruence topology} on $M$. We put
$$ C(M)= \displaystyle{\bigcap_{\frak{q} \neq
\{0\}}}\hat{M}(\frak{q}),
$$
where $\hat{M}(\frak{q})$ is the usual profinite completion of
$M(\frak{q})$ with respect to all its finite index subgroups. We
call $C(M)$ the {\it congruence kernel} of $M$. Then there is an
exact sequence of the above type involving $C(M)$,
$\hat{M}(\frak{q})$ and the completion of $M(\frak{q})$ with
respect to the congruence topology.

\medskip

We may assume that $\Gamma$ and, hence all its subgroups, act on
the Bruhat-Tits tree $T$ associated with $G$ {\it without
inversion}. As usual the vertex and edge sets of a graph $X$ will
be denoted by $V(T)$ and $E(T)$, respectively. Given a subgroup
$H$ of $\Gamma$ and $w\in V(T)\cup E(T)$, we denote by $ H_w$ the
stabiliser of $w$ in $H$ Since $\Gamma$ is discrete it follows
that $H_w$ is always finite.

\medskip

 We deal with the cocompact and non-uniform cases separately.

\section{ Cocompact arithmetic lattices}

For each positive integer $s$, let $F_s$ denote the free group of
rank $s$.

\begin{lem} If $\Gamma$ is cocompact, then, for all but
finitely many $\frak{q}$,
$$
\Gamma(\frak{q}) \cong F_r,
$$
where $r = r(\frak{q}) \geq 2$. \noindent Moreover $r (\frak q)$
is unbounded in the following sense.\\If $r(\frak{q}) \geq 2$  and
$$
\frak{q} = \frak{q}_1 \gneqq \frak{q}_2 \gneqq \frak{q}_3 \cdots
$$
is an infinite properly descending chain of $\mathcal{O}$-ideals,
then
$$
r(\frak{q}_i) \rightarrow \infty, \: as\;\: i \rightarrow \infty.
$$
\end{lem}

\noindent {\bf Proof.}  It is well-known that the quotient graph
$\Gamma \backslash T $  is  finite. Let $v_1, \cdots, v_t$ denote
the vertices (in $V(T)$) of a {\it lift} $j: \Gamma \backslash T
\rightarrow T$. We put
$$
\Gamma_i = \Gamma_{v_i} \qquad\, (1 \leq i \leq t).
$$
It is clear that, for all but finitely many $\frak{q}$,
$$
\Gamma(\frak{q}) \cap \Gamma_i = \{I_n\} \qquad \,(1 \leq i \leq
t),
$$
since each $\Gamma_i$ is finite. For such a $\frak{q}$ all the
stabilizers in $\Gamma(\frak{q})$ of the vertices of $T$ are
trivial, since $\Gamma(\frak{q})$ is normal in $\Gamma$.  Further
$|\Gamma: \Gamma(\frak{q})|$ is finite and so $\Gamma(\frak{q})
\backslash T$ is finite. It follows that
$$
\Gamma(\frak{q}) \cong F_r,
$$
for some $r$; see [S2, Theorem 4, p.~27].  By Lemma~1 it is clear
that $r \geq 2$. If $r(\frak{q}) \geq 2$ and
$$
\frak{q} = \frak{q}_1 \gvertneqq \frak{q}_2 \gvertneqq \frak{q}_3
\cdots $$
then by the well-known Schreier formula,
$$
r(\frak{q}_i)-1 = |\Gamma(\frak{q}):
\Gamma(\frak{q}_i)|(r(\frak{q})-1).
$$
The result follows since $|\Gamma(\frak{q}): \Gamma(\frak{q}_i)|
\rightarrow \infty$, as $i \rightarrow \infty$. $\hfill\Box$\\

\begin{thm} If $\Gamma$ is cocompact, then
$$
C(\Gamma) \cong \hat{F}_{\omega}.
$$
\end{thm}
\noindent {\bf Proof.}  Fix any $\frak{q}$ for which Lemma 2.1
holds. Let $C =C(\Gamma)$. Then, by the exact sequence after Lemma
1.3,
$$
\hat{F}_r/C \cong \overline{\Gamma}(\frak{q}).
$$
Now $|\mathbf{G}(\mathcal{O}): \Gamma(\frak{q})|$ is finite and so
(by Lemma 1.2)  $\overline{\Gamma}(\frak{q})$ embeds as an {\it
open} subgroup of $\mathbf{G}(\hat{\mathcal{O}})$ and hence
contains an open subgroup $O$ of $\mathbf{G}(\hat{\mathcal{O}})$
of type (*).

\noindent Since $\Gamma$ is cocompact, $\Gamma(\frak{q})$ is
finitely generated.  It follows that
$\overline{\mathbf{G}}(\mathcal{O}), \overline{\Gamma}(\frak{q})$
and $O$ are all {\em finitely generated\/} profinite groups.
Consequently the group $O$ does not "satisfy Schreier's formula''.
(See [RZ, Lemma 8.4.5, p.~320].) Hence
$\overline{\Gamma}(\frak{q})$ does not satisfy Schreier's formula,
since $|\overline{\Gamma}(\frak{q}): O|$ is finite. The result
follows from [RZ, Corollary 8.4.4, p.~320]. $\hfill\Box$

\section{ Non-uniform arithmetic lattices: discrete results}

Here we assume that $G/\Gamma$ is {\it not} compact, in which case
$k$ is a function field. We put $\char k = p$. It is well-known
that an element $X$ of $\Gamma$ has finite order if and only if $X
\in \Gamma_v$, for some $v \in V(T)$. In order to describe the
structure of $\Gamma \backslash T$ we make the following

\medskip

\noindent {\bf Definition.} Let $R$ be a {\it ray} in $\Gamma
\backslash T$, i.e. an infinite path without backtracking and let
$j:R \rightarrow T$ be a {\it lift}. Let $V(j(R))=\{v_1,
v_2,\cdots\}$. We say that $j$ is {\it stabilizer ascending}, if $
\Gamma_{v_i} \leq \Gamma_{v_{i+1}}$ for $i \geq 1$,  and set
$$
\Gamma(R)\;(=\Gamma(R, j)) := \langle\Gamma_v\,|\,\, v\in
V(j(R))\rangle.
$$
\noindent Using results of Raghunathan [R], Lubotzky [L2, Theorem
6.1] has determined the structure of $\Gamma\backslash T$.  This
extends an earlier result of Serre [S, Theorem 9, p.~106] for the
special case $\mathbf{G} =\mathbf{SL}_2,\quad \Gamma = {\rm
SL}_2(\mathcal{O})$ and $S = \{v\}$. Baumgartner [Ba] has provided
a more detailed and extended version of Lubotzky's proof.

\begin{thm} With the above notation,
$$
\Gamma \backslash T = Y \cup R_1 \cup \cdots \cup R_m,
$$
where $Y$ is a finite subgraph and $R_1, \cdots, R_m$ are rays. In
addition,
\begin{itemize}
\item[(a)] $\mathrm{card}\,\{V(Y) \cap V(R_i)\} = 1 \qquad (1 \leq
i \leq m)$, \item[(b)] $E(Y) \cap E(R_i) = \emptyset \qquad (1
\leq i \leq m)$, \item[(c)] $R_i \cap R_\ell = \emptyset \qquad (i
\neq \ell)$.
\end{itemize}
There exists a lift $j: \Gamma \backslash T \rightarrow T$ such
that $ j: R_i \rightarrow T$ is stabilizer ascending for $1 \leq i
\leq m$.
\end{thm}

\smallskip

\begin{lem} With the above notation, the group $\Gamma(R_i)$ is contained
in $\mathbf{P}_i(k_v)$, where $\mathbf{P}_i$ is a minimal
parabolic $k_v$-subgroup of $\mathbf{G}$, where $1 \leq i \leq m.$
\end{lem}

\noindent {\bf Proof.} The group $\Gamma(R_i)$ stabilizes the end
of $T$ corresponding to $j(R_i)$. It is well-known from the
standard theory of Bruhat-Tits that the stabilizer of an end in
$G$ is of the form $\mathbf{P}_i(k_v)$. $\hfill\Box$

\medskip

We now restrict our attention to principal congruence subgroups.

\begin{lem} Let $\mathfrak{q}$ be a proper
$\mathcal{O}$-ideal. Then every element of finite order of
$\Gamma(\mathfrak{q})$ is unipotent of $p$-power order.
\end{lem}

\noindent{\bf Proof.}  Let $k_0$ be the (full) field of constants
of (the function field) $k$. Let $g \in \Gamma(\mathfrak{q})$ have
finite order and let $\chi_g(t)$ denote its characteristic
polynomial over $k$. Then

$$
\chi_g(t) \equiv (t-1)^n \ (\mod \mathfrak{q}).
$$
Now each zero of $\chi_g(t)$ is a root of unity and so each
coefficient of $\chi_g(t)$ lies in the algebraic closure of $k_0$
in $k$, which is $k_0$ itself. Since $k_0 \leq \mathcal{O}$ it
follows that $ \chi_g(t) = (t-1)^n. $ $\hfill\Box$

\smallskip

\begin{lem} With the above notation, for each proper
$\mathcal{O}$-ideal $\frak{q}$, let
$$
\Gamma(\frak{q}) \cap \Gamma(R_i)= \Theta_i(\frak{q})
$$
and let $\mathbf{U}_i$ be the unipotent radical of $\mathbf{P}_i$,
where $1 \leq i \leq m$ . Then
\begin{itemize}
\item[(i)] $\Theta_i(\mathfrak{q})$ is a subgroup of finite index in $\mathbf{U}_i(\mathcal{O})$; \item[(ii)]
$\Theta_i(\frak{q})$ is nilpotent of class at most $2$ and
 is
generated by elements of p-power order..
\end{itemize}
\end{lem}

\noindent{\bf Proof.}  Since $\Theta_i(\frak{q})$ consists of
elements of finite order in $\Gamma(\mathfrak{q})$ it consists of
unipotent matrices by Lemma 3.3. Part (i) follows. (Recall that
$\Gamma$ is an {\it arithmetic} lattice.) For part (ii) we note
that $\mathbf{G}$ has $k_v$-rank one and so $\mathbf{U}_i$ is
nilpotent of class at most $2$, by [BT2, 4.7 Proposition].
$\hfill\Box$

\medskip

As we shall see some (but not all) such $\mathbf{U}_i$ are in fact
abelian.

\begin{thm} For all but finitely many $\mathfrak{q}$,
$$
\Gamma(\mathfrak{q}) \cong F_r * \Lambda(\mathfrak{q}),
$$
where $\Lambda(\mathfrak{q})$ is a free product of finitely many
groups, each of which is a conjugate (in $\Gamma$) of some
$\Theta_i(\frak{q})$. (Then $\Lambda(\frak{q})$ is
generated by nilpotent groups of class at most $2$, each consisting of elements of p-power order.)\\[.25cm]
In addition,
$$
r = r(\mathfrak{q}) = {\mathrm
rk}_{\mathbb{Z}}(\Gamma(\mathfrak{q}))=\mathrm{dim}_\mathbb{Q}
H^1(\Gamma(\mathfrak q),\mathbb{Q}),
$$
the (finite) free abelian rank of $\Gamma(\mathfrak{q})$.
\end{thm}

\noindent {\bf Proof.}  By the fundamental theorem of the theory
of groups acting on trees [S2, Theorem 13, p.~55] $\Gamma$ is the
fundamental group of the graph of groups given by the lift $j:
\Gamma \backslash T \rightarrow T$ as described in Theorem 3.1.
For all but finitely many $\mathfrak{q}$,
$$
\Gamma(\mathfrak{q}) \cap \Gamma_v = \{I_n\},
$$
for all $v \in V(j(Y))$. We fix such a $\mathfrak{q}$. Recall that
$\Gamma(\mathfrak{q})$ is a {\it normal} subgroup of finite index
in $\Gamma$. From standard results on the decomposition of a
normal subgroup of a fundamental group of a graph of groups,
$\Gamma(\mathfrak{q})$ is a free product of a free group $F_r$ and
a finite number of subgroups, each of which is a conjugate of
$\Gamma(\mathfrak{q})\cap\Gamma(R_i)$, for some $i$. The rest
follows from Lemma 3.4. $\hfill\Box$

\medskip

\noindent For the case $\mathbf{G}=\mathbf{SL}_2,\,
S=\{v\}\;\mathrm{and}\; \Gamma={\rm SL}_2(\mathcal{O})$, Theorem
3.5 is already known [Mas2, Theorem 2.5].

\begin{cor} Let $U(\mathfrak{q})$ denote the (normal)
subgroup of $\Gamma(\mathfrak{q})$ generated by its unipotent
matrices. Then, for all but finitely many $\mathfrak{q}$,
$$
\Gamma(\mathfrak{q})/U(\mathfrak{q}) \cong F_r,
$$
where $r = r(\mathfrak{q}) = {\mathrm
rk}_{\mathbb{Z}}(\Gamma(\mathfrak{q}))$.
\end{cor}

\noindent {\bf Proof.} We fix an ideal $\mathfrak{q}$ for which
Theorem 3.5 holds. Let $\Lambda(\mathfrak{q})^*$ denote the normal
subgroup of $\Gamma(\mathfrak{q})$ generated by
$\Lambda(\mathfrak{q})$. Now every unipotent element of
$\Gamma(\mathfrak{q})$ is of finite order and so lies in a
conjugate of some $\Theta_i(\mathfrak{q})$, by Theorem 3.5. It
follows that
$\Lambda(\mathfrak{q})^*=U(\mathfrak{q})$.$\hfill\Box$

\medskip

We now show that $r(\frak{q})$ is not bounded.

\begin{lem} With the above notation, for infinitely many
$\mathfrak{q}$ we have $$r(\mathfrak{q}) \geq 2.$$ If
$r(\mathfrak{q}') \geq 2$ and $ \mathfrak{q}' = \mathfrak{q}_1
\gvertneqq\mathfrak{q}_2\gvertneqq
\mathfrak{q}_3\gvertneqq\cdots\, $ is an infinite properly
descending chain of $\mathcal{O}$-ideals, then
$$
r(\mathfrak{q}_i) \rightarrow \infty, \; \mbox{ as }\; i
\rightarrow \infty.
$$
\end{lem}

\medskip

\noindent{\bf Proof.}  We note that, if $ \Gamma(\frak{q}) = F_s *
H, $ where $H$ is a subgroup of $\Gamma(\mathfrak{q})$, then $
r(\mathfrak{q}) \geq s. $ By Theorem 3.1 together with [S, Theorem
13, p.~55] it follows that $ \Gamma = A *_{W} B, $ where
\begin{itemize}
\item[(i)] $B = \Gamma(R)$, for some ray $R$ and a lift $j: R
\rightarrow T$; \item[(ii)] $W = \Gamma_v$, for some $v \in V(T)$.
\end{itemize}
Now $B$ is infinite (since $\Gamma$ is {\it non-uniform}) and $W$
is finite. If $A=W$, then $\Gamma(\mathfrak{q})$ is nilpotent by
Lemma 3.4, for any proper $\mathfrak{q}$.  This contradicts Lemma
1.1.  We conclude that $W\neq A$.

\medskip

 It is well-known that, for any
$\mathfrak{q}$,
$$
r(\mathfrak{q}) \geq 1 + |\Gamma: W\cdot\Gamma(\mathfrak{q})| -
|\Gamma: A\cdot\Gamma(\mathfrak{q})| - |\Gamma: B\cdot
\Gamma(\mathfrak{q})|.
$$
We now restrict our attention to the (all but finitely many)
$\frak{q}$ for which
$
W \cap \Gamma(\mathfrak{q}) = \{I_n\}.
$
Among these are infinitely many $\mathfrak{q}'$ for which
$$
|A \cdot \Gamma(\mathfrak{q}'):\Gamma(\mathfrak{q}')| > |W \cdot
\Gamma(\mathfrak{q}'):\Gamma (\mathfrak{q}')| \; \mbox{ and }\; |B
\cdot \Gamma(\mathfrak{q}'):\Gamma(\mathfrak{q}')|> 2 |W \cdot
\Gamma(\mathfrak{q}'):\Gamma(\mathfrak{q}')|.
$$
It follows that $ r(\mathfrak{q}') \geq 2. $ For the second part,
it is clear that
$$
r(\mathfrak{q}_{i+1}) \geq r(\mathfrak{q}_i) \geq 2 \qquad (i \geq
1).
$$
Fix $i$.  Then, by Theorem 3.5, $ \Gamma(\mathfrak{q}_i) = F_{r'}
* H, $ say, where $r' = r(\mathfrak{q}_i)$.  For any $t > i$, it
follows from the Kurosh subgroup theorem and the Schreier formula
that $ r(\mathfrak{q}_t) > r', $ unless $\Gamma(\mathfrak{q}_t)
\cap F_{r'} = F_{r'}$ and $\Gamma(\mathfrak{q}_i) =
\Gamma(\mathfrak{q}_t) \cdot F_{r'}$. We choose $t$ so that
$\Gamma(\mathfrak{q}_i) \neq \Gamma(\mathfrak{q}_t)$. $\hfill\Box$

\medskip

Lemma 3.7 is already known for the case
$\mathbf{G}=\mathbf{SL}_2,\, S=\{v\}$ and $ \Gamma={\rm
SL}_2(\mathcal{O})$. See the proof of [Mas1, Theorem 3.6].

\medskip

 Before providing a complete description of $C(\Gamma)$ for the
non-uniform case we first establish a special property of
unipotent groups in rank one algebraic groups.\\
\section{The congruence kernel of a unipotent group}
\noindent We assume that $\mathbf{G}$, $k$, $\mathcal{O}$ and
$k_v$ are as above. Let $\mathbb K$ be an algebraically closed
field containing $k_v$. In view of Theorem 2.2 we will assume from
now on that $k$ is a {\it function field}, with
$\mathrm{char}\;k=p$. (Although a number of results in this
section also hold for number fields.)  Throughout $\mathbf{P}$
denotes a minimal $k_v$-parabolic subgroup of $\mathbf{G}$ and
$\mathbf{U}$ denotes its unipotent radical (also defined over
$k_v$). Let $U=\mathbf{U}(\mathcal{O})$. Now by [BT, 4.7
Proposition] it follows that the congruence kernel
$$
C(U)= \displaystyle{\bigcap_{\frak{q} \neq
\{0\}}\hat{U}(\frak{q})}
$$
is nilpotent of class at most $2$. The principal aim of this
Section is to prove that $C(U)$ is in fact {\it abelian}.\\

\noindent We note that since $\mathbf{G}$ is $k$-isotropic it has
$k$-rank one. Making use of [PRag], it follows from Tits
Classification [T] that $\mathbf{G}$ is isomorphic to one of a
(finite) number of types. In Tits notation [T] (adapted) we
conclude that $\mathbf{G}$ is isomorphic to one of the following:
\begin{itemize}
\item[(a)] Inner type $\mathbf{A}$;\item[(b)] Outer type $\mathbf{A}_{2d+1}$ which becomes inner over $k_v$; \item[(c)] Outer types
$\mathbf{A}_2,\; \mathbf{A}_3$; \item[(d)] Types $\mathbf{C}_2,\;
\mathbf{C}_3$; \item[(e)] Types $\mathbf{D}_3,\; \mathbf{D}_4,\;
\mathbf{D}_5$.
\end{itemize}

\noindent Now if $\mathbf{G}$ is an inner form of type
$\mathbf{A}$ then ${\mathbf G}(k_v)=\mathbf{SL}_2(D)$ where $D$ is
a central simple division algebra over $k_v$.  In this case it is
known that $\mathbf{U}$ is abelian. This is also true when
$\mathbf{G}$ is of type $\mathbf{C}_2$. (See [PRag, 1.1, 1.3] for
more details.) For case (b) the groups can be realised as
two-dimensional special unitary groups over a division algebra $D$
of degree $d$ with centre $K$, where $K$ is a separable quadratic
extension of $k$. (The description of the groups involves an
involution (of the second kind) which is defined on $D$.) Now the
place $v$ of $k$ splits over $K$ and so, over $k_v$, $\mathbf{G}$
is of inner type $\mathbf{A}_{2d+1}$. Here then $\mathbf{U}$ is
also abelian.  (See [PR, p.352].) For the purposes of this Section
therefore we need not consider these cases any further. For outer
forms of type ${\bf A}_2$ and ${\bf A}_3\cong {\bf D}_3$ the Tits
indices are

\input epsf

\begin{center}
$\epsfbox{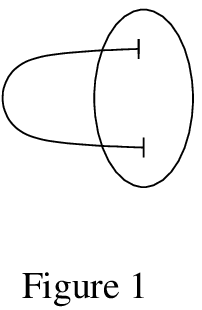}\qquad\qquad\qquad {\epsfbox{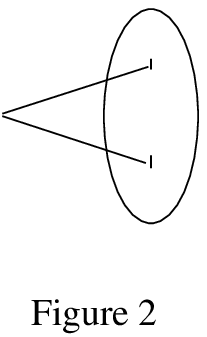}}$
\end{center}
while for type ${\bf C}_3$ it has the following
form:\\

\begin{center}

\epsfbox{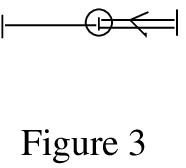}

\end{center}

\noindent Finally for types ${\bf D}_4$ and ${\bf D}_5$ the
indices are
\begin{center}
$\qquad \epsfbox{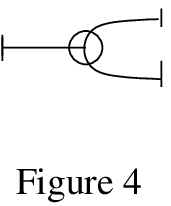}\qquad\qquad\qquad\
\epsfbox{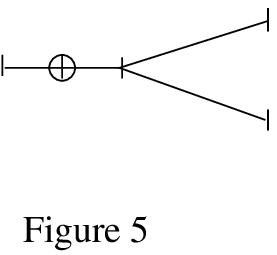}$
\end{center}

We now recall some generalities on reductive algebraic $k$-groups
which will be useful later on. Let ${\mathcal G}={\mathbf
G}(\mathbb K)$ and ${\mathfrak g}={\rm Lie}\,{\mathcal G}$, the
Lie algebra of the algebraic group ${\mathcal G}$. Let
$X_*({\mathcal G})$ denote the set of all cocharacters of
$\mathcal G$, i.e. the set of all rational homomorphisms from the
multiplicative group $\mathbb K^\times$ to $\mathcal G$. Note that
for any $\lambda\in X_*({\mathcal G})$ the group $\lambda(\mathbb
K^\times)$ is a $1$-dimensional torus in $\mathcal G$.

Given $\phi \in X_*(\mathcal G)$ and $g\in \mathcal G$ we say that
the limit $\lim_{t\rightarrow 0} \phi(t)g\phi(t^{-1})$ {\it
exists} if the morphism from ${\mathbb K}^\times$ to ${\mathcal
G}$ sending $t\in{\mathbb K}^\times$ to
$\phi(t)g\phi(t^{-1})\in\mathcal G$ extends to a morphism from
$\mathbb K$ to $\mathcal G$. Let
\begin{eqnarray*} {\mathcal P}(\phi)&:=&\{g\in {\mathcal G}\,|
\,\,\lim_{t\rightarrow 0}
\phi(t)g\phi(t^{-1})\,\mbox{ exists}\,\}\\
 {\mathcal U}(\phi)&:=&\{g\in {\mathcal G}\,|\,\,\lim_{t\rightarrow 0}
\phi(t)g\phi(t^{-1})=1\}.
\end{eqnarray*}
\noindent It is well-known that ${\mathcal P}(\phi)$ is a
parabolic subgroup of $\mathcal G$ and ${\mathcal U}(\phi)$ is the
unipotent radical of ${\mathcal P}(\phi)$. Moreover, if the
morphism $\phi$ is defined over $k$, then both ${\mathcal
P}(\phi)$ and ${\mathcal U}(\phi)$ are $k$-defined subgroups of
$\mathcal G$; see [Sp, I, 4.3.4 and II, 3.3.1].

\medskip

Crucial for our purposes is the following surprising result. It
ensures that the structure of {\it any} $C(U)$ can deduced from a
detailed description of {\it one particular} $\mathbf{U}$. (This
result in fact holds for any global field.)

\begin{thm}
If  $U$ is nonabelian, then $\mathbf{U}$ is defined over $k$.
\end{thm}
\noindent {\bf Proof.} Let ${\mathcal P}={\mathbf P}(\mathbb K)$
and ${\mathcal U}=\mathbf{U}(\mathbb K)$. Obviously, $\mathcal P$
is a parabolic subgroup of $\mathcal G$ and ${\mathcal
U}=R_{u}({\mathcal P})$, the unipotent radical of $\mathcal P$.
Choose a maximal torus $\mathcal T$ of $\mathcal G$ contained in
$\mathcal P$ and let $\Phi$ denote the root system of $\mathbf G$
relative to $\mathcal T$. Denote by $X({\mathcal T})$  the lattice
of rational characters of $\mathcal T$, and let $\Delta$ be a
basis of simple roots in $\Phi$. Adopt Bourbaki's numbering of
simple roots and denote by $\widetilde{\alpha}$ the highest root
of $\Phi$ with respect to $\Delta$.

Let $\alpha^\vee$ denote the coroot corresponding to
$\alpha\in\Phi$, an element in $X_*({\mathcal T})\subset
X_*({\mathcal G})$. Recall $\alpha^\vee({\mathbb K}^\times)$ is a
$1$-dimensional torus in $\mathcal T$. As usual, we let ${\mathcal
U}_\alpha=\{x_\alpha(t)\,|\,\,t\in\mathbb K\}$ denote the root
subgroup of $\mathcal G$ corresponding to $\alpha$; see [St, \S3].
Given $x\in{\mathcal G}$ we denote by $Z_{\mathcal G}(x)$ the
centraliser of $x$ in $\mathcal G$.

\smallskip

\noindent {\bf Case} $\mathbf{1}$. We first suppose that $\mathbf
G$ is not of type ${\mathbf C}_3$. The above discussion then shows
that $\mathbf G$ is of type ${\mathbf A}_2$, ${\mathbf A}_3$,
${\mathbf D_4}$ or ${\mathbf D}_5$. A quick look at the Tits
indices displayed above reveals that $\mathcal P$ is $\mathcal
G$-conjugate to the normaliser in $\mathcal G$ of the
$1$-parameter unipotent subgroup ${\mathcal
U}_{\widetilde{\alpha}}$. From this it follows that in our present
case the derived subgroup of $\mathcal U$ has dimension $1$ as an
algebraic group and coincides with the centre $\mathcal Z$ of
$\mathcal U$. Moreover, ${\mathcal Z}$ is ${\mathcal G}$-conjugate
to ${\mathcal U}_{\widetilde{\alpha}}$.

\medskip

By our assumption, the derived subgroup $[U,U]$ contains an
element $u\ne 1$. Then $$u\in [U,U]\subset [{\mathbf
U}(k_v),{\mathbf U}(k_v)]\subset [{\mathcal U},{\mathcal
U}]=\mathcal Z.$$ Since the subgroup ${\mathcal Z}$ is $\mathcal
T$-invariant, the preceding remark implies that there is a long
root $\beta\in\Phi$ such that ${\mathcal U}={\mathcal U}_\beta$.
Then $u=x_\beta(a)$ for some {\it nonzero} $a\in {\mathbb K}$. We
claim that the centraliser $Z_{\mathcal G}(u)$ is defined over
$k$. To prove this claim it suffices to verify that the orbit
morphism $g\longmapsto  g u g^{-1}$ of $\mathcal G$ is separable;
see [Sp, II, 2.1.4]. The latter amounts to showing that the Lie
algebra of $Z_{\mathcal G}(u)$ coincides with ${\mathfrak
g}^u:=\{X\in {\mathfrak g}\,|\,\, ({\rm Ad}\, u)(X)=X\}$.

\medskip

After adjusting $\Delta$, possibly, we can assume that
$\beta=\widetilde{\alpha}$. For each $\alpha\in\Phi$ we choose a
nonzero vector $X_\alpha$ in ${\mathfrak g}_\alpha={\rm
Lie}\,{\mathcal U}_\alpha$ and let ${\mathfrak t}$ denote the Lie
algebra of $\mathcal T$. Denote by ${\mathfrak g}'$ the $\mathbb
K$-span of all $X_\gamma$ with
$\gamma\not\in\{\pm\widetilde{\alpha}\}$ and set ${\mathfrak
g}(\widetilde{\alpha}):={\mathfrak
g}_{-\widetilde{\alpha}}\oplus{\mathfrak t}\oplus{\mathfrak
g}_{\widetilde{\alpha}}$. Clearly, ${\mathfrak g}={\mathfrak
g}'\oplus{\mathfrak g}(\widetilde{\alpha})$. Using [St, \S3] it is
easy to observe that both ${\mathfrak g}'$ and ${\mathfrak
g}(\widetilde{\alpha})$ are $({\rm Ad}\,u)$-stable and one can
choose $X_{\widetilde{\alpha}}$ such that
$$({\rm
Ad}\,u)(X_\gamma)=X_\gamma+a[X_{\widetilde{\alpha}},X_\gamma]
\qquad\ (\forall \,\,X_\gamma\in{\mathfrak g}').$$ Since
$\widetilde{\alpha}$ is long, standard properties of root systems
and Chevalley bases imply that if $\gamma\in\Phi$ is such that
$\gamma\ne -\widetilde{\alpha}$ and
$\gamma+\widetilde{\alpha}\in\Phi$, then $[X_{\widetilde{\alpha}},
X_\gamma]=\lambda_\gamma X_{\widetilde{\alpha}+\gamma}$ for some
{\it nonzero} $\lambda_\gamma\in{\mathbb K}$; see [St, Theorem~1].
From this it follows that ${\mathfrak g}^u\cap{\mathfrak g}'$
coincides with the $\mathbb K$-span of all $X_\gamma$ such that
$\gamma\not\in\{\pm\widetilde{\alpha}\}$ and
$\widetilde{\alpha}+\gamma\not\in\Phi$. On the other hand, the
commutator relations in [St, Lemma~15] imply that each such
$X_\gamma$ belongs to ${\rm Lie}\,Z_{\mathcal G}(u)$. Therefore,
${\mathfrak g}^u\cap{\mathfrak g}'\subset{\rm Lie}\,Z_{\mathcal
G}(u).$

\medskip

The differential ${\rm d}\widetilde{\alpha}$ is a linear function
on ${\mathfrak t}$. Since ${\mathbf G}$ is simply connected, the
equality ${\rm d}\widetilde{\alpha}=0$ holds if and only if
$\widetilde{\alpha}=p\mu$ for some $\mu\in X(\mathcal T)$. The
latter holds if and only if $p=2$ and ${\mathbf G}$ is of type
${\mathbf A}_1$ or ${\mathbf C}_n$. Thus, in the present case we
have that ${\rm d}\widetilde{\alpha}\ne 0$. As
$$({\rm Ad}\,u)(h)=h-a({\rm
d}\widetilde{\alpha})(h)X_{\widetilde{\alpha}}\qquad\, (\forall\,
h\in\mathfrak t),$$ this implies that ${\mathfrak
g}^u\cap{\mathfrak g}(\widetilde{\alpha})={\mathfrak
g}_{\widetilde{\alpha}}\oplus \ker {\rm d}\widetilde{\alpha}$. But
then ${\mathfrak g}^u\cap{\mathfrak
g}(\widetilde{\alpha})\subset{\rm Lie}\,Z_{\mathcal G}(u),$
forcing ${\mathfrak g}^u\subseteq{\rm Lie}\,Z_{\mathcal G}(u)$.
Since ${\rm Lie}\,Z_{\mathcal G}(x)\subseteq{\mathfrak g}^x$ for
all $x\in\mathcal G$, we now derive that the group $Z_{\mathcal
G}(u)$ is defined over $k$. Hence the connected component (of the
identity of) $Z_{\mathcal G}(u)^\circ$ is defined over $k$, too;
see [Sp, II, 2.1.1].

\medskip

Let $\mathcal C$ denote the connected component of the centraliser
$Z_{\mathcal G}({\mathcal U}_{\widetilde{\alpha}})$. The argument
above shows that ${\rm Lie}\, {\mathcal C}={\rm Lie}\,Z_{\mathcal
G}(u)$. Since ${\mathcal C}\subseteq Z_{\mathcal G}(u)^\circ$, we
must have the equality $Z_{\mathcal G}(u)^\circ=\mathcal C$. Then
$\mathcal C$ is a $k$-group, hence contains a maximal torus
defined over $k$, say ${\mathcal T}'$. As
$\ker\widetilde{\alpha}\subset Z_{\mathcal G}(u)$, the torus
${\mathcal T}'$ has dimension $l-1$, where $l={\rm rk}\,\mathbf
G$. Let $\mathcal H$ denote the centraliser of ${\mathcal T}'$ in
$\mathcal G$. By construction, $\mathcal H$ is a connected
reductive $k$-group of semisimple rank $1$ containing ${\mathcal
U}_{\widetilde{\alpha}}$. Since $\mathcal G$ is simply connected,
so is the derived subgroup of ${\mathcal H}$; see [SS, II,
Theorem~5.8].  As ${\mathcal U}_{\widetilde{\alpha}}$ is
unipotent, it lies in $[{\mathcal H},{\mathcal H}]$. As $1\ne u\in
{\mathbf G}(k)\cap [{\mathcal H},{\mathcal H}]$, the group
$[{\mathcal H},{\mathcal H}]$ is $k$-isotropic. The classification
of simply connected $k$-groups of type ${\bf A}_1$ now shows that
$[{\mathcal H},{\mathcal H}]\cong {\rm SL}_2({\mathbb K})$ as
algebraic $k$-groups. As a consequence, $u$ belongs to a
$k$-defined Borel subgroup of $[{\mathcal H},{\mathcal H}]$; call
it $\mathcal B$. Since $u$ commutes with ${\mathcal
U}_{\widetilde{\alpha}}$, it must be that ${\mathcal
U}_{\widetilde{\alpha}}=R_u({\mathcal B})$.

 \medskip

Let $\mathcal S$ be a $k$-defined maximal torus of $\mathcal B$.
Since $[{\mathcal H},{\mathcal H}]$ is $k$-isomorphic to ${\rm
SL}_2({\mathbb K})$, there exists a $k$-defined cocharacter
$\mu\colon\,{\mathbb K}^\times \rightarrow [{\mathcal H},{\mathcal
H}]$ such that $${\mathcal S}=\mu({\mathbb K}^\times),\qquad
\mu(t)x_{\widetilde{\alpha}}(t')\,\mu(t)^{-1}=x_{\widetilde{\alpha}}(t^2t')
\quad\,\,(\forall\,t, t'\in\mathbb K).$$ Then $\mathcal S$ is
$k$-split in $\mathcal G$, and hence it is a {\it maximal}
$k_v$-split torus of $\mathcal G$ (recall that $\mathbf G$ has
$k_v$-rank $1$). Since $\mathcal S$ normalises ${\mathcal
U}_{\widetilde{\alpha}}$, it lies in ${\mathcal P}$. As $\mathcal
P$ is defined over $k_v$, there exists a $k_v$-defined cocharacter
$\nu\colon\,{\mathbb K}^\times\rightarrow \mathcal P$ such that
${\mathcal P}={\mathcal P}(\nu)$; see [Sp, II, 5.2.1]. Since
$\mu({\mathbb K}^\times)$ and $\nu({\mathbb K}^\times)$ are
maximal $k_v$-split tori in $\mathcal P$, they are conjugate by an
element of $\mathcal U$; see [Sp, II, Theorem~5.2.3 (iv)]. In
conjunction with the earlier remarks this yields that $r\nu={\rm
Int}\, x\circ s\mu$ for some $x\in\mathcal U$ and some positive
integers $r$ and $s$. But then $${\mathcal P}(\mu)={\mathcal
P}(s\mu)={\mathcal P}({\rm Int}\, x\circ s\mu)={\mathcal
P}(r\nu)={\mathcal P}.$$ Since $\mu$ is defined over $k$, so are
$\mathbf P$ and $\mathbf U$; see [Sp, II, 3.1.1.].

\medskip

\noindent {\bf Case} $\mathbf{2}$. Next suppose that $\mathbf G$
is of type ${\mathbf C}_3$ and $p\ne 2$. As before, we denote by
$\nu$ a $k_v$-defined cocharacter in $X_*(\mathcal P)$ such that
${\mathcal P}={\mathcal P}(\nu)$.  Let ${\mathcal G}^{\rm uni}$
and ${\mathfrak g}^{\rm nil}$ denote the unipotent variety of
$\mathcal G$ and the nilpotent variety of $\mathfrak g,$
respectively. These are affine varieties defined over $k$. Since
$\mathcal G$ is simply connected and $p$ is a good prime for
$\Phi$, the Bardsley--Richardson projection associated with a
semisimple $k$-representation of $\mathcal G$ induces a
$k$-defined, $\mathcal G$-equivariant isomorphism of varieties
$$\eta\,\colon\,\,{\mathcal G}^{\rm
uni}\stackrel{\sim}{\longrightarrow}{\mathfrak g}^{\rm nil}$$ such
that $\eta({\mathcal U})={\rm Lie}\,\mathcal U$; see [McN, 8.5]
for more detail. Set $X:=\eta(u)$, a $k$-rational nilpotent
element of $\mathfrak g$. Since $X$ is an unstable vector of the
$({\rm Ad}\,\mathcal G)$-module $\mathfrak g$, associated with $X$
is a nonempty subset $\tilde{\Lambda}_X\subset X_*({\mathcal G})$
consisting of the so-called {\it optimal} cocharacters for $X$;
see [P, 2.2] for more detail. Since in the present case the orbit
map $g\longmapsto ({\rm Ad}\,g)(X)$ of $\mathcal G$ is separable
at $X$, by [SS, I, \S 5] for example, it follows from the main
results of [McN] that $\tilde{\Lambda}_X$ contains a $k$-defined
cocharacter $\lambda$ such that $({\rm Ad}\,\lambda(t))(X)=t^2 X$
for all $t\in{\mathbb K}^\times$. Since $u\in [{\mathcal
U},{\mathcal U}]$, it is immediate from Figure~3 and the
definition of $\eta$ that $({\rm Ad}\,\nu)(t))=t^{2m}X$ for some
positive integer $m$. But then $\lambda(t)^m\nu(t)^{-1}\in
Z_{\mathcal G}(X)$ for all $t\in\mathbb K^\times$, where
$Z_{\mathcal G}(X)=\{g\in{\mathcal G}\,|\,\,({\rm Ad}\,g))(X)=X\}$
is the centraliser of $X$ in $\mathcal G$.

\medskip
Since $\lambda$ gives an optimal torus for $X$, the instability
parabolic subgroup ${\mathcal P}(\lambda)$ contains $Z_{\mathcal
G}(X)$; see [P, 2.2] for example. Since $\lambda$ is defined over
$k$, so is ${\mathcal P}(\lambda)$; see [Sp, II, 3.1.1.]. As
$\lambda({\mathbb K}^\times)\subset {\mathcal P}(\lambda)$, the
preceding remark yields $\nu({\mathbb K}^\times)\subset {\mathcal
P}(\lambda)$.

 Since $\nu({\mathbb K}^\times)$ and $\lambda({\mathbb K}^\times)$ are maximal
 $k_v$-split tori in
${\mathcal P}(\lambda)$, they are conjugate in ${\mathcal
P}(\lambda)$; see [Sp, II, 5.2.3]. It follows that there exists
$x\in{\mathcal P}(\lambda)$ such that $r\nu={\rm Int}\,x\circ
s\lambda$ for some positive integers $r$ and $s$. But then
$r\nu\in\tilde{\Lambda}_X$; see [P, 2.2] for example. As a result,
$${\mathcal P}={\mathcal P}(\nu)={\mathcal P}
(r\nu)={\mathcal P}(\lambda).$$ Since $\lambda$ is defined over
$k$, so are $\mathbf P$ and $\mathbf U$, see [Sp, II, 3.3.1].

\medskip

\noindent {\bf Case} $\mathbf{3}$. Finally, suppose that $\mathbf
G$ is of type ${\mathbf C}_3$ and $p=2$. In this case we cannot
argue as in {\bf Case} $\mathbf{2}$ because $p=2$ is  bad for
$\Phi$. We shall argue as in {\bf Case} $\mathbf{1}$ instead. Let
$\beta_0=\alpha_1+2\alpha_2+\alpha_3=\varepsilon_1+\varepsilon_2$,
the highest {\it short} root in $\Phi$, and
\begin{eqnarray*}
\Gamma_0&:=&\{\pm\alpha_1,\pm\alpha_3\},\\
\Gamma_1&:=&\{\alpha_2,\alpha_1+\alpha_2,
\alpha_2+\alpha_3,\alpha_1+\alpha_2+\alpha_3\},\\
\Gamma_2&:=&\{2\alpha_2+\alpha_3,
\beta_0,2\alpha_1+2\alpha_2+\alpha_3\}.
\end{eqnarray*}
According to Figure~3, it can be assumed that $\mathcal U$ is
generated by the unipotent root subgroups ${\mathcal U}_\gamma$
with $\gamma\in\Gamma_1\cup\,\Gamma_2$. Moreover, $\langle\,
{\mathcal U}_\gamma\,|\,\,\gamma\in\Gamma_2\rangle$ is a central
normal subgroup of $\mathcal U$ containing the derived subgroup of
$\mathcal U$. Furthermore, $\mathcal P$ is generated by $\mathcal
T$, $\mathcal U$, and $\langle\, {\mathcal
U}_\gamma\,|\,\,\gamma\in\Gamma_0\rangle.$

\medskip

Since $p=2$, combining the above description of ${\mathcal U}$
with Steinberg's relations [St, Lemma 15] shows that $[{\mathcal
U},{\mathcal U}]={\mathcal U}_{\beta_0}$ and $\mathcal P$
coincides with the normaliser of ${\mathcal U}_{\beta_0}$ in
$\mathcal G$. It follows that $[U,U]$ contains an element
$u=x_{\beta_0}(a)$ for some nonzero $a\in\mathbb K$. Consequently,
$$\dim\,{\rm Lie}\,Z_{\mathcal G}(u)=\dim\,Z_{\mathcal
G}(u)=\dim\,{\mathcal P}-1.$$ We adopt the notation of $\mathfrak
t$, $X_\gamma$, ${\mathfrak g}_\gamma$, and ${\mathfrak g}^u$
introduced in {\bf Case~2}. For $i\in\{\pm 1,\pm 2\}$, we denote
by ${\mathfrak g}_i$ the $\mathbb K$-span of all $X_\gamma$ with
$\gamma\in\pm\Gamma_i$, and let ${\mathfrak g}_0$ be the ${\mathbb
K}$-span of $\mathfrak t$ and all $X_\gamma$ with
$\gamma\in\Gamma_0$. Then ${\rm Lie}\,{\mathcal P}=\bigoplus_{i\ge
0}\,{\mathfrak g}_i$. The decomposition
$${\mathfrak g}={\mathfrak g}_{-2}\oplus{\mathfrak
g}_{-1}\oplus{\mathfrak g}_0\oplus{\mathfrak g}_1\oplus{\mathfrak
g}_2$$ gives $\mathfrak g$ a graded Lie algebra structure. In view
of [St, \S3] we have
$$\big({\rm Ad}\,u-{\rm Id}\big)({\mathfrak g}_k)\subseteq
\,\textstyle{\bigoplus}_{i\ge 2}\,{\mathfrak g}_{k+i}\qquad\ \,
(\forall\,\, k\ge -2).$$ Take $x\in{\mathfrak g}^u$ and write
$x=\sum\,x_i$ with $x_i\in{\mathfrak g}_i$. Combining [St, \S3]
with the preceding remark it is straightforward to see that
$$0\equiv \big({\rm Ad}\,u-{\rm
Id}\big)(x_{-2})\equiv\, a[X_{\beta_0},x_{-2}]\ \,\, \big(\mod\,
\textstyle{\bigoplus}_{i\ge 1}\,{\mathfrak g}_i\big).$$ On the
other hand, standard properties of Chevalley bases (and the fact
that $\mathbf G$ is simply connected) ensure that ${\rm
ad}\,X_{\beta_0}$ is injective on ${\mathfrak g}_{-2}$. Therefore,
$x_{-2}=0$. Arguing similarly we obtain $x_{-1}=0$. As a result,
${\mathfrak g}^u\subseteq {\rm Lie}\,{\mathcal P}$.

\medskip

Similar to {\bf Case} $\mathbf{1}$ we observe that the
differential ${\rm d}\beta_0$ is a nonzero linear function on
${\mathfrak t}$. As $({\rm Ad}\,u)(h)=h-a({\rm
d}\beta_0)(h)X_{\beta_0}$ for all $h\in\mathfrak t$, this implies
that ${\mathfrak g}^u$ is a {\it proper} Lie subalgebra of ${\rm
Lie}\,\mathcal P$. But then $\dim\,{\mathfrak g}^u=\dim\,{\rm
Lie}\,Z_{\mathcal G}(u)$, forcing ${\rm Lie}\,Z_{\mathcal
G}(u)={\mathfrak g}^u$. Hence $Z_{\mathcal G}(u)$ is defined over
$k$. Then so is the connected component of $Z_{\mathcal G}(u)$;
see [Sp, II, 2.1.1].

\medskip

We now denote by $\mathcal C$ denote the connected component of
the centraliser $Z_{\mathcal G}({\mathcal U}_{\beta_0})$. The
above argument shows that ${\rm Lie}\, {\mathcal C}={\rm
Lie}\,Z_{\mathcal G}(u)$. Then $Z_{\mathcal G}(u)^\circ=\mathcal
C$, so that $\mathcal C$ is a $k$-group. We let ${\mathcal T}'$ be
a maximal $k$-defined torus in ${\mathcal C}$ and denote by
$\mathcal H$ the centraliser of ${\mathcal T}'$ in $\mathcal G$.
At this point we can repeat verbatim our argument in {\bf Case~1}
to conclude that there is a $k$-defined cocharacter
$\mu\colon\,{\mathbb K}^\times \rightarrow [{\mathcal H},{\mathcal
H}]$ such that $\mu({\mathbb K}^\times)$ normalises ${\mathcal
U}_{\beta_0}$. Our earlier remarks then yield ${\mu({\mathbb
K}^\times)\subset \mathcal P}$. As in {\bf Case~1} this implies
that both $\mathbf P$ and $\mathbf U$ are defined over $k$. This
completes the proof. $\hfill\Box$

\medskip

\noindent {\bf Remark}. Let $L/F$ be a field extension and let
$\bf G$ be an absolutely simple, simply connected algebraic
$F$-group. Suppose that ${\rm char}\,F$ is either zero or a very
good prime for $\bf G$ (the list of very good primes is well-known
and can be found in [McN, 2.1] for example). Suppose further that
$\bf G$ has $L$-rank $1$ and let $\bf P$ be a minimal parabolic
subgroup of $\bf G$ defined over $L$. Let $\bf U$ be the unipotent
radical of $\bf P$ and suppose that  $$[{\bf U}(L),{\bf U}(L)]\cap
{\mathbf G}(F)\ne \{1\}.$$ Then it follows from the argument used
in {\bf Case} $\mathbf{2}$ of the proof of Theorem~4.1 that ${\bf
U}$ is defined over $F$. (One should also keep in mind that $[{\bf
U},[{\bf U},{\bf U}]]=\{1\}$, which one can see by analyzing the
list of Tits indices in [Sp, pp.~81-83].) Our proof of Theorem~4.1
suggests that that this might even be true without any
restrictions on the characteristic of $F$.

\begin{lem}
Let $\mathbf{P}_i$ be a minimal parabolic $k$-subgroup of
$\mathbf{G}$ with unipotent radical $\mathbf{U}_i$ and let
$\mathbf{U}_i(\mathcal{O})= U_i$, where $i=1,2$. Then
$$
C(U_1) \cong C(U_2).
$$
\end{lem}
\noindent {\bf Proof.} By standard Borel-Tits theory
$\mathbf{P}_1$, $\mathbf{P}_2$ (and hence $\mathbf{U}_1$,
$\mathbf{U}_2$) are conjugate over $k$. The result follows from
[Mar, Lemma 3.1.1, p.~60]. $\hfill\Box$

\medskip

\noindent Our next result is especially important. We recall from
[BT, 4.7 Proposition] that $[U,U]$ is a central subgroup of $U$.
\begin{lem} Let $Z=Z(\mathcal{O})$ be a (possibly trivial) central subgroup of
$U$, containing the commutator subgroup $[U,U]$, such that $U/Z$
is a countably infinite elementary abelian $p$-group. Suppose
further that, if $N$ is any subgroup of finite index in $U$, then
$$ Z(\frak{q}) \leq N,
$$
for some non-zero $\mathcal{O}$-ideal $\frak{q}$. \newline Then
$C(U)$ is isomorphic to the direct product of $2^{{\aleph}_0}$
copies of $\mathbb{Z}/p\mathbb{Z}$.
\end{lem}
\noindent {\bf Proof.}  Let $C=C(U)$ and
$\mathfrak{c}=2^{{\aleph}_0}$. Since any vector space of countably
infinite dimension has $\mathfrak{c}$ hyperplanes, the hypotheses
ensure that $U$ has $\mathfrak{c}$ finite index subgroups. On the
other hand $\mathcal{O}$ has only countably many ideals and so $U$
has ${\aleph}_0$ congruence subgroups. It follows that
$$
\mathrm{card}(C)=2^{\mathfrak{c}}.
$$
The hypotheses also ensure that
$$
C\cap \bar{Z} = \{1\},
$$ where $\bar{Z}$ denotes the closure of $Z$ in $\hat{U}$.
It follows that $C$ embeds in
$$
\hat{U}/\bar{Z} \cong \hat{V},
$$
where $V=U/Z$. The result follows. $\hfill\Box$\\

\noindent Note that Lemma 4.3 applies to the case where $U$ is a
countably infinite elementary abelian $p$-group. For the remainder
of this section we say that any $U$ with a central subgroup $Z$
satisfying the hypotheses in the statement of Lemma 4.3 {\it has
property} $\mathfrak{P}$.  We now proceed to prove that this lemma
applies to all $C(U)$ on a case-by-case basis.

\medskip

\noindent {\bf Case} $\mathbf{1}$: {\bf Outer types}
$\mathbf{A}_2,\mathbf{A}_3$

\medskip

\noindent Let $K$ be a (Galois) quadratic extension of $k$, and
let $\sigma$ be the generator of the Galois group of $K/k$. Let
$f$ be the $\sigma$-hermitian, non-degenerate form in $n+1$
variables over $K$ determined by the matrix
$$
F=\pmatrix{0 & 0 & 1 \cr 0 & F_0 & 0 \cr 1 & 0 & 0},
$$
where (i) $F_0=1$, when $n=2$, and (ii) $F_0$ is a
$\sigma$-hermitian, anisotropic $2\times 2$ matrix over $K$, when
$n=3$. As usual, for any matrix $M$ over $K$, we put $M^* =
(M^{\sigma})^{\mathrm{tr}}$. For $n=2,3$ we define
$$
\mathbf{SU}(K,f):=\{ X \in {\rm SL}_{n+1}(K)\,|\,\, X^*FX=F\}.
$$
Clearly we can represent this group in ${\rm SL}_{2n+2}(k)$ by
means of any $2$-dimensional representation of $K$ over $k$. The
following is an immediate consequence of [T].

\begin{thm} Let $\mathbf{G}$ be of outer type $\mathbf{A}_n$ where
$n=2,3$.
Then there exist $K, f$ of above type such that
$$
G(=\mathbf{G}(k)) \cong \,\mathbf{SU}(K,f).
$$

\end{thm}

\noindent We now denote by $\mathbf{UT}(K,f)$ the set of all upper
unitriangular matrices in ${\rm SL}_{2n+2}(k)$ contained in
$\mathbf{SU}(K,f)$.

\begin{lem}
If $\mathbf{G}$ is of outer type $\mathbf{A}_n$ where $n=2,3$,
then there exists a minimal parabolic $k$-subgroup $\mathbf{P}_0$
of $\mathbf{G}$ with unipotent radical $\mathbf{U}_0$, such that
$$
\mathbf{U}_0(k) \cong \mathbf{UT}(K,f). $$
\end{lem}
\noindent {\bf Proof.} First, let us consider $\mathbf G$ of outer
type ${\mathbf A}_2$. Let $K/k$ and $\sigma$ be as before, and let
$A$ be any commutative algebra over $k$. Then $\sigma$ extends
uniquely to an $A$-linear involution on the $K$-algebra
$A\otimes_k K$. Let $G(A)\,=\,\{g \in {\mathbf SL}_3(A\otimes_k
K)\,|\,\, g^{\ast}Fg = F \}$, where $g^\ast=(g^\sigma)^{\rm tr}$
and $$F = \pmatrix{0 & 0 & 1 \cr 0 & 1 & 0 \cr 1 & 0 & 0}.$$ It
follows from the Tits classification that $G(A)$ is the group of
$A$-rational points of a simple algebraic $k$-group $k$-isomorphic
to $\mathbf G$. Thus we may assume without loss of generality that
${\mathbf G}({\mathbb K})=G({\mathbb K})$.

\medskip

Identify $\mathbb K$ with ${\mathbb K}\otimes_k k\subset {\mathbb
K}\otimes_k K$, and define  $\nu\in X_*\big({\mathbf
SL}_3({\mathbb K}\otimes_k K)\big)$ by setting
$$\nu(t) := \mathrm{diag}(t,\, 1,\,t^{-1})
\qquad (\forall\,t \in {\mathbb K}^{\times}).$$ Put
$S:=\nu({\mathbb K}^\times)$. As $S\subset {\bf G}(\mathbb K)$, we
have that $\nu\in X_*(G(\mathbb K))$. The above description of $G$
yields that the morphism $\nu\colon\,{\mathbb K}^\times\rightarrow
G(\mathbb K)$ is defined over $k$.

\medskip

Direct computation shows that the parabolic subgroup of ${\mathbf
SL}_3({\mathbb K}\otimes_k K)$ associated with $\nu$ is nothing
but the group of all upper triangular matrices in ${\mathbf
SL}_3({\mathbb K}\otimes_k K)$. In other words $(P(\nu))(\mathbb
K)$ is nothing but the group of all upper triangular matrices in
$G(\mathbb K)$. As a consequence, the unipotent radical of
$(P(\nu))(\mathbb K)$ coincides with the group of all upper
unitriangular matrices in $G(\mathbb K)$. More precisely, for
$\alpha, \beta, \gamma \in {\mathbb K}\otimes_k K$ define
$$T(\alpha, \beta, \gamma)\, :=\, \pmatrix{1 & \alpha & \beta \cr 0 & 1
& \gamma \cr 0 & 0 & 1}.$$

\noindent Then $(U(\nu))({\mathbb
K})\,=\,\{T(\alpha,\beta,\gamma)\,|\,\, \gamma = -
\alpha^{\sigma}, \beta + \beta^{\sigma} = - \alpha
\alpha^{\sigma}\}.$ Since $G$ has $k$-rank $1$, the group
$(U(\nu))(\mathbb K)$ must be equal to the unipotent radical of a
minimal $k$-parabolic subgroup of $G(\mathbb K)$.

\medskip

We consider outer type ${\bf A}_3$ now. In this case also, $K/k$
and $\sigma$ are as before, and for any commutative algebra $A$
over $k$, $\sigma$ extends uniquely to an $A$-linear involution on
the $K$-algebra $K\otimes_k A$. The group $G(A)\,=\,\{g \in
{\mathbf SL}_4(K\otimes_k A)\,|\,\, g^{\ast}Fg = F \}$, where
$g^\ast=(g^\sigma)^{\rm tr}$ and
$$F = \pmatrix{0 & 0 & 0 & 1 \cr
0 & a & b & 0 \cr 0 & b^{\sigma} & d & 0 \cr 1 & 0 & 0 & 0}.$$

\noindent From the Tits classification, we have that $G(A)$ is the
group of $A$-rational points of a simple algebraic $k$-group
$k$-isomorphic to $\mathbf G$. Thus we may assume without loss of
generality that ${\mathbf G}({\mathbb K})=G({\mathbb K})$.

Identifying $\mathbb K$ with ${\mathbb K}\otimes_k k\subset
{\mathbb K}\otimes_k K$, we get a cocharacter $\nu\in
X_*\big(\mathbf{SL}_4({\mathbb K}\otimes_k K)\big)$ by setting
$$\nu(t) := \mathrm{diag}(t,\, 1, 1, \,t^{-1})
\qquad (\forall\,t \in {\mathbb K}^{\times}).$$ Put
$S:=\nu({\mathbb K}^\times)$. As $S\subset G(\mathbb K)$, we have
that $\nu\in X_*(G(\mathbb K))$. The above description of $G$
yields that the morphism $\nu\colon\,{\mathbb K}^\times\rightarrow
G(\mathbb K)$ is defined over $k$. Exactly, as in the case of
${\bf A}_2$, an easy computation shows that the unipotent radical
of the (minimal) $k$-parabolic subgroup associated to $\nu$
consists of the upper unitriangular matrices in $G(\mathbb K)$.
 $\hfill\Box$\\

 For $n=2,3$ we denote the $(n+1) \times (n+1)$ matrix
$$
\pmatrix{1 & \alpha & \beta \cr 0 & 1 & \gamma \cr 0 & 0 & 1}
$$ by $T(\alpha,\beta, \gamma)$, where $\alpha$ and $\beta^{\mathrm{tr}}$
are $1 \times (n-1)$. We note that
$T(\alpha_1,*,\gamma_1)T(\alpha_2,*,\gamma_2)=T(\alpha_1+\alpha_2,*,\gamma_1+\gamma_2)$.
\begin{lem}
$$
\mathbf{UT}(K,f) = \{ T(\alpha,\beta,\gamma)\in
\mathbf{SU}(K,f)\,|\,\, \alpha = -\gamma^*F_0, \;\beta
+\beta^{\sigma} = -\gamma^*F_0\gamma\}.
$$
\end{lem}
\noindent {\bf Proof.} We note that any $2\times 2$ unipotent
matrix over $k$ representing an element of $K$ is the identity. In
addition the {\it only}  upper unitriangular matrix $Y$ over $K$
such that $Y^*F_0Y=F_0$ is the
identity, since $F_0$ is anisotropic. $\hfill\Box$\\

\noindent The following is readily verified.
\begin{lem}
Suppose that $T(*,\beta_i,\gamma_i) \in \mathbf{UT}(K,f)$, where
$i=1,2$. Then
\begin{itemize}
\item[(a)] $T(*,x^2\beta_1,x\gamma_1) \in \mathbf{UT}(K,f)$, for
all $x\in k$,
\item[(b)]$[T(*,\beta_1,\gamma_1),T(*,\beta_2,\gamma_2)]=
T(0,\lambda-\lambda^{\sigma},0)$,
where $\lambda= \gamma_1^*F_0\gamma_2$.
\end{itemize}
\end{lem}

\noindent The $k$-subspace of $K$
$$ V=\{s-s^{\sigma}\,|\,\, s\in K\}
$$
has $k$-dimension $1$. In choosing a pair of $2\times2$ matrices
(with entries in $\mathcal{O}$) as a $k$-basis for $K$, we ensure
that one of them, $\mathbf{v}$, say, spans $V$. With the notation
of Lemma 4.5 we put $UT=\mathbf{U}_0(\mathcal{O})$.
\begin{lem}  $UT$ has property
$\mathfrak{P}$.
\end{lem}
\noindent {\bf Proof.} There exist $T(*,*,\gamma_i)\in UT$, where
$i=1,2$, such that $\gamma_1^*A\gamma_2 - \gamma_2^*A\gamma_1 \neq
0$. Now let $N$ be any finite index normal subgroup of $UT$. Then
by Lemma 4.7(a) we may assume that $T(*,*,\gamma_1)\in N$.
\noindent It is easily verified from Lemmas 4.6, 4.7 (a) that
$$
Z (=Z(\mathcal{O}))=\{ T(0,y\mathbf{v},0)\,|\,\, y\in
\mathcal{O}\}
$$
is a (non-trivial) central subgroup of $UT$, containing $[UT,UT]$.
Now $N\cap Z$ then contains $[T(*,*,\gamma_1),T(*,*,y\gamma_2)]$,
for all $y\in \mathcal{O}$. It follows that $Z(\frak{q}) \leq N$,
for some non-zero (principal) $\mathcal{O}$-ideal, $\mathfrak{q}$.
It is clear form the above that $UT/Z$ is (infinite) elementary $p$-abelian.$\hfill\Box$\\

\noindent {\bf Case} $\mathbf{2}$: {\bf Type} $\mathbf{C}_3$\\

\noindent Let $D$ be a quaternion division algebra over $k$ and
let $\sigma$ be an involution of $D$ of the first kind (i.e. an
anti-homomorphism of $D$ of order $2$ which fixes $k$.). Suppose
that $D^{\sigma}$, the $k$-subspace of $D$ containing all elements
of $D$ fixed by $\sigma$, has $k$-dimension $3$. Let $h$ be the
$\sigma$-skewhermitian, non-degenerate form in $3$ variables over
$D$ determined by the matrix $$ H=\pmatrix{0 & 0 & 1 \cr 0 & d & 0
\cr -1 & 0 & 0},
$$
where $d^{\sigma}=-d \neq 0.$ We define
$$
\mathbf{SU}(D,h)=\{X \in {\rm SL}_3(D)\,|\,\, X^*HX=H \}.
$$
Clearly we can represent this group in ${\rm SL}_{12}(k)$ by means
of any $4$-dimensional representation of $D$ over $k$. The
following is an immediate consequence of [T].
\begin{thm}
Let $\mathbf{G}$ be of type $\mathbf{C}_3$. Then there exist $D,h$
of the above type such that
$$
G(=\mathbf{G}(k))\cong \mathbf{SU}(D,h).
$$
\end{thm}

\noindent As above we consider the subgroup $\mathbf{UT}(D,h)$ of
all upper unitriangular matrices in ${\mathrm SL}_{12}(k)$
contained in $\mathbf{SU}(D,h)$.

\begin{lem}
 There exists a minimal parabolic $k$-subgroup $\mathbf{P}_0$ of
 $\mathbf{G}$ with unipotent radical $\mathbf{U}_0$ such that
 $$
 \mathbf{U}_0(k)\cong \mathbf{UT}(D,h).
 $$
 \end{lem}
 \noindent {\bf Proof.} The proof will be similar to that of Lemma~4.5.
Here, $D$ is a quaternion division algebra with an involution
$\sigma$ of the first kind, and $G$ is the special unitary group
of a non-degenerate $\sigma$-skew-hermitian form on a
$3$-dimensional (right) $D$-vector space. The form can be
represented by the matrix $$\pmatrix{0 & 0 & 1 \cr 0 & d & 0 \cr
-1 & 0 & 0},\quad\ d\in D^\times,\ \ d^{\sigma} = -d.$$ We get a
rational homomorphism $\nu\colon\,{\mathbb K}^\times\rightarrow
G({\mathbb K})={\mathbf SL}_3({\mathbb K}\otimes_k D)$ by setting
$$\nu(t) := \mathrm{diag}(t,\, 1,\,t^{-1})
\qquad (\forall\,t \in {\mathbb K}^{\times}).$$ It is defined over
$k$ and $S:=\nu({\mathbb K}^\times)$ is a maximal $k$-split torus
of $G$. The rest of the proof is as before.$\hfill\Box$

\medskip

 Continuing with the above notation we use $T(\alpha,\beta,\gamma)$
 to denote this time
 a $3\times 3$ upper unitriangular matrix over $D$,
 where $\alpha, \beta, \gamma \in D$.
 \begin{lem}
 $$
 \mathbf{UT}(D,h) = \{T(\alpha,\beta,\gamma)\in \mathbf{SU}(D,h):
 \alpha=\gamma^{\sigma}d,\; \beta -\beta^{\sigma}=
 \gamma^{\sigma}d\gamma\}.
 $$
 \end{lem}
 \noindent {\bf Proof.} We note that the {\it only} unipotent matrix
 over $k$ representing an element of $D$ is the identity.
 $\hfill\Box$

\medskip

 Lemma 4.8 has the following equivalent.

 \begin{lem}
 Suppose that $T(*,\beta_i,\gamma_i) \in \mathbf{UT}(D,h)$, where
 $i=1,2$. Then
 \begin{itemize}
 \item[(a)] $T(*,x^2\beta^i,x\gamma_i)\in \mathbf{UT}(D,h)$, for all
 $x\in k$,
 \item[(b)]
 $[T(*,\beta_1,\gamma_1),T(*,\beta_2,\gamma_2)]=T(0,\lambda+\lambda^{\sigma},0)$,
 where $\lambda= \gamma_1^{\sigma}d\gamma_2$.
 \end{itemize}
 \end{lem}
 \noindent As we see later for our purposes this case is
 essentially identical to that of type $\mathbf{D}_3$,
 when $\mathrm{char}\;k=2$, by Lemma 4.11. For now therefore we assume
 that $\mathrm{char}\;k \neq 2$.
 \noindent The $k$-subspace of $D$
 $$
 \{ x\in D\,|\,\, x^{\sigma}=-x\}
 $$
 has $k$-dimension $1$. We may choose four $4\times4$ matrices over
 $k$, $\mathbf{v}_i$, where $i=1,2,3,4$, as a $k$-basis for $D$,
 with $\mathbf{v}_i^{\sigma}=\mathbf{v}_i$, when $i=1,2,3$, and
 $\mathbf{v}_4=d$. We may assume that all the entries of these
 matrices lie in $\mathcal{O}$. By considering $(d^3)^{\sigma}$ it
 is clear that $d^2=\mu$, for some (non-zero) $\mu\in\mathcal{O}$.
 the following is very easily verified.
 \begin{lem}
 When $i=1,2,3$
 $$
 [T(*,*,r_i\mathbf{v}_i),T(*,*,s_i\mathbf{v}_4)]=T(0,2r_{i}s_i\mu\mathbf{v}_i,0),$$
 for all $r_i,s_i\in k$.
 \end{lem}

\noindent As before we put $UT=\mathbf{U}_0(\mathcal{O})$ in the
notation of Lemma 4.10.
\begin{lem}Suppose that $\mathrm{char}\;k\neq 2$.
Then $UT$ has property $\mathfrak{P}$.
\end{lem}
\noindent {\bf Proof.} We note that by Lemma 4.11 the element
$T(*,*,2r\mathbf{v}_i)\in UT$, for all $r\in\mathcal{O}$, where
$i=1,2,3,4$. Let
$$
Z(=Z(\mathcal{O}))=\{T(0,\beta,0)\in UT: \beta^{\sigma}=\beta\}.
$$
Then from the above $Z$ is a central subgroup of $UT$, containing
$[UT,UT]$. Let $N$ be a normal subgroup of finite index in $UT$.
>From the above $T(*,*,r_i\mathbf{v}_i)\in N$, for some non-zero
$r_i\in\mathcal{O}$. Let $r_0=r_1r_2r_3.$ Then
$$
T(0,2s_1r_0\mu\mathbf{v}_1+2s_2r_0\mu\mathbf{v}_2+2s_3r_0\mu\mathbf{v}_3)\in
N\cap Z,
$$
for all $s_1,s_2,s_3\in \mathcal{O}$. It follows that
$Z(\mathfrak{q})\leq N$, for some non-zero (principal)
$\mathcal{O}$-ideal, $\mathfrak{q}$. It is clear from the above that $UT/Z$ is an (infinite)
 elementary abelian $p$-group.$\hfill\Box$\\

\noindent {\bf Case} $\mathbf{3}:$ {\bf Types}
$\mathbf{D}_3,\mathbf{D}_4,\mathbf{D}_5$\\

\noindent Let $D,\; \sigma$ be as above. Let $q$ be a
$\sigma$-quadratic, non-degenerate form in $n$ variables over $D$
and let $q'$ be its associated $\sigma$-hermitian form, where
$n=3,4,5$. Suppose further that $q$ has Witt index $1$ over $k$.
(When $\mathrm{char}\;k=2$ it is assumed also that $q$ is {\it
non-defective}.)
\begin{thm}
Let $\mathbf{G}$ be of type $\mathbf{D}_n$, where $n=3,4,5$. Then
there exists $q$ of the above type and a central $k$-isogeny
$$
\pi:\mathbf{G}\rightarrow \mathbf{SO}(q).
$$
In addition, if $\mathbf{U}$ is the unipotent radical of a minimal
parabolic $k$-subgroup of $\mathbf{G}$, then $\pi(\mathbf{U})$ is
the unipotent radical of a minimal parabolic $k$-subgroup of
$\mathbf{SO}(q)$ which is $k$-isomorphic to $\mathbf{U}$.
\end{thm}

\noindent {\bf Proof.} Follows from [T] and [BT1, Propositions
2.20,2.24]. $\hfill\Box$\\

\noindent We now represent $q'$ by means of the $n\times n$ matrix
over $D$
$$
L=\pmatrix{0 & 0 & 1 \cr 0 & Q & 0 \cr 1 & 0 & 0},
$$
where $Q$ is an $(n-2)\times (n-2)$ anisotropic,
$\sigma$-hermitian matrix. Then the $k$-rational points of
$\mathbf{SO}(q)$ are given by
$$
 \mathbf{SU}(D,q') =\{X\in {\rm SL}_n(D)\,|\,\,X^*LX=L\}.
$$
\noindent As before we can use any $4\times4$ representation of
$D$ over $k$ to obtain a $4n\times4n$ representation of
$\mathbf{SU}(D,q')$ over $k$. \noindent We let $\mathbf{UT}(D,q')$
denote the subgroup of all upper unitriangular matrices in ${\rm
SL}_{4n}(k)$  contained in $\mathbf{SU}(D,q')$. Adapting a
previous notation we put
$$
T(\alpha,\beta,\gamma)= \pmatrix{1 & \alpha & \beta \cr 0 & 1 &
\gamma \cr 0 & 0 & 1},
$$
where $\alpha, \beta^{\mathrm{tr}}$ are matrices of type $1\times
(n-2)$ over $D$ ($n=3,4,5$).
\begin{lem} There exists a minimal parabolic $k$-subgroup of
$\mathbf{G}$ with unipotent radical $\mathbf{U}_0$, such that
$$
\mathbf{U}_0(k)\cong\mathbf{UT}(D,q').
$$
\end{lem}
\noindent {\bf Proof.} We shall replace $G$ by (and work with) the
image of $G$ under the central $k$-isogeny in 4.15. Thus, we have
a quaternion division algebra $D$, an involution $\sigma$ of the
first kind, and an $n \times n$ matrix ($n=3,4,5$)
$$L = \pmatrix{0 & 0 & 1 \cr 0 & Q & 0 \cr 1 & 0 & 0}$$
where $Q$ is an $(n-2) \times (n-2)$ matrix which represents a
$\sigma$-hermitian, anisotropic form.\\
We are working with the subgroup of ${ \mathrm SL}_n(D)$ which
preserves $L$. In this case the rational homomorphism is:
$$\nu\colon\,{\mathbb K}^\times\rightarrow\, G({\mathbb K})=\,\mathbf{SL}_n({\mathbb
K}\otimes_k D), \quad\, t \mapsto \mathrm{diag}(t,\, 1, \ldots,
1,\,t^{-1}) \qquad (\forall\,t \in {\mathbb K}^{\times}).$$ The
size of the matrix is $3,4$ or $5$, according as we are in
$D_3,D_4$ or $D_5$. In all cases, the proof is similar.
$\hfill\Box$

\begin{lem}
$$
\mathbf{UT}(D,q') =\, \{T(\alpha,\beta,\gamma)\in
\mathbf{SU}(D,q')\,|\,\, \alpha = -\gamma^*Q,\;
\beta+\beta^{\sigma}=-\gamma^*Q\gamma\}.
$$
\end{lem}
\noindent {\bf Proof.} As before the only unipotent matrix over
$k$ representing an element of $D$ is the identity. In addition
the only upper triangular unipotent matrix $W$ over $D$, such that
$W^*QW=Q$ is again the identity.$\hfill\Box$

\medskip

Lemmas 4.7 and 4.12 have the following equivalent.

\begin{lem}
 Suppose that $T(*,\beta_i,\gamma_i)\in \mathbf{UT}(D,q')$, where
 $i=1,2$. Then
 \begin{itemize}
 \item[(a)] $T(*,x^2\beta_i,x\gamma_i)\in\mathbf{UT}(D,q')$, for
 all $x\in k$,
 \item[(b)]
 $[T(*,\beta_1,\gamma_1),T(*,\beta_2,\gamma_2)]=T(0,\lambda-\lambda^{\sigma},0)$,
 where $\lambda=\gamma_1^*Q\gamma_2$.
\end{itemize}
\end{lem}

\noindent The hypotheses on $D$ ensure that the $k$-subspace of
$D$
$$
\{ d-d^{\sigma}\,|\,\, d \in D\}
$$
has $k$-dimension $1$. We can therefore choose a $k$-basis of $D$,
consisting of four $4\times4$ matrices, with entries in
$\mathcal{O}$, one of which spans this subspace. Let
$UT=\mathbf{U}_0(\mathcal{O})$, where $\mathbf{U}_0$ is as defined
in Lemma 4.16. From the above, in a way very similar to Lemma 4.8
we can prove the following.
\begin{lem}
$UT$ has property $\mathfrak{P}$.
\end{lem}
\noindent We note that since Lemma 4.19 includes type
$\mathbf{D}_3$, Lemma 4.14 also holds (for type $\mathbf{C}_3$)
when $\mathrm{char}\;k=2$. We now come to the main conclusion of
this section.

\begin{thm}
Let $\mathbf{U}$ be the unipotent radical of a minimal parabolic
$k_v$-subgroup of $\mathbf{G}$ and let
$U=\mathbf{U}(\mathcal{O})$. Then the congruence kernel $C(U)$ is
isomorphic to the direct product of $2^{{\aleph}_0}$ copies of
$\mathbb{Z}/p\mathbb{Z}$.
\end{thm}
\noindent {\bf Proof.} There are two possibilities. If $U$ is
abelian then, from [T], $\mathbf{G}$ is either inner type
$\mathbf{A}$ or type $\mathbf{C}_2$. From [PRag, 1.1, 1.3] and
standard Borel-Tits theory it follows that $U$ is an elementary
abelian $p$-group. We can now apply Lemma 4.3. \newline
Alternatively $\mathbf{U}$ is defined over $k$ by Theorem 4.1. The
result follows from Lemmas 4.2, 4.3, 4.8, 4.14 and 4.19.
$\hfill\Box$

\section{ Non-uniform arithmetic lattices: profinite results}

Continuing from the previous section we assume that $k$ is a
function field with $\mathrm{char}\; k=p$. Let $A$ and $B$ be
profinite groups. We will denote by
$$
A\amalg B
$$
the {\it free profinite product of A and B}. See [RZ, p.~361].

\medskip

 \noindent Let $\hat{F}_s$ denote {\it the free profinite
group of (finite) rank s}, where $s\geq 1$.

\begin{lem} With the above notation, for all but finitely many
$\mathfrak{q}$,
$$
\hat{\Gamma}(\mathfrak{q})\cong\hat{F}_r\amalg\hat{\Lambda}(\mathfrak{q}),
$$
where
\begin{itemize}
\item[(a)] $\hat{\Lambda}(\mathfrak{q})$ is a free profinite
product of nilpotent pro-$p$ groups, each of which is of the type
$\hat{\Theta}(\mathfrak{q})$, where
$$\Theta(\mathfrak{q})=\Gamma\cap\mathbf{U}(\mathfrak{q}),
$$
for some unipotent radical $\mathbf{U}$ of a minimal
$k_v$-parabolic subgroup of $\mathbf{G}$. (In which case
$\hat{\Theta}(\mathfrak{q})$ is nilpotent of class at most $2$ and
is generated by torsion elements of $p$-power order.);
\item[(b)] the normal subgroup of $\hat{\Gamma}(\mathfrak{q})$
generated by $\hat{\Lambda}(\mathfrak{q})$ is
$\hat{U}(\mathfrak{q})$ ;
\item[(c)]$r=r(\mathfrak{q})$ is not bounded.
\end{itemize}
Moreover,
$$
\hat{\Gamma}(\mathfrak{q})/\hat{U}(\mathfrak{q})\cong\hat{F}_r.
$$
\end{lem}
\noindent{\bf Proof.} Follows from Theorem 3.5 and Lemma
3.7.$\hfill\Box$

\medskip

A {\it projective} group is, by definition, a closed subgroup of a
free profinite group.
\begin{lem}
Let $N$ be a normal, closed, non-open subgroup of
$\hat{\Gamma}(\mathfrak{q})$. Then, for all but finitely many
$\mathfrak{q}$,
$$
N\cong P\amalg N(\mathfrak{q}),
$$
where
\begin{itemize}
\item[(a)] $N(\mathfrak{q})$ is a closed subgroup of $\hat{U}(\mathfrak{q})$
and a free profinite product of
nilpotent pro-$p$ groups, each of class at most $2$ and each
generated by torsion elements of $p$-power order;
\item[(b)] $P$ is a projective group, all of whose proper, open subgroups
are isomorphic to $\hat{F}_{\omega}$.
\end{itemize}
\end{lem}
\noindent {\bf Proof.} This follows from a result of the fourth
author [Za1, Theorem 2.1]. (See also [Za1, Theorem 4.1, Lemma
4.2].)$\hfill\Box$

\medskip

An immediate consequence of Lemma 5.2 and Lemma 1.3 is the
following.
\begin{lem}
With the above notation,
$$
C(\Gamma)\cong P\amalg N(\Gamma),
$$
where
\begin{itemize}
\item[(a)]  $N(\Gamma)$ is
a closed subgroup of all  $\hat{U}(\mathfrak{q})$ and a free
profinite product of elementary abelian pro-$p$ groups;
\item[(b)] $P$ is a projective group, all of whose proper, open
subgroups are isomorphic to $\hat{F}_{\omega}$.
\end{itemize}
\end{lem}
\noindent{\bf Proof.} We apply Lemma 5.1 and the proof of Lemma
5.2 to the case $N=C(\Gamma)$. Then $C(\Gamma)$ is the free
profinite product of $P$, as above, and (in the notation of Lemma
5.1) groups of the type $C(\Gamma)\cap\hat{\Theta}(\mathfrak{q})$.
By Lemmas 1.3 and 3.4 it follows that
$$
C(\Gamma) \cap\hat{\Theta}(\mathfrak{q})=
\displaystyle{\bigcap_{\frak{q}' \neq
\{0\}}\hat{\Gamma}(\mathfrak{q}')\cap\hat{\Theta}(\frak{q})}=\displaystyle{\bigcap_{\{0\}\neq\frak{q}'\leq\frak{q}}\hat{\Theta}(\frak{q}')}\leq
C(U).
$$
The result follows from Theorem 4.20. $\hfill\Box$\\

\noindent {\bf Terminology.} If we can replace $P$ with
$\hat{F}_{\omega}$ in Lemma 5.3, we will say that the {\it principal
result} holds.
\begin{lem}
Let $A$ and $B$ be profinite groups and let $M$ be a normal,
closed subgroup of
$$
A\amalg B.
$$
Then $M\cap A$ is a factor in the free profinite decomposition of
$M$.
\end{lem}
\noindent {\bf Proof.} Follows from [Za1, Theorem
2.1].$\hfill\Box$
\begin{lem}
Let $P$ be as in Lemma 5.3 and $F$ be isomorphic to
$\hat{F}_{\omega}$. Then
$$
P\amalg F\cong \hat{F}_{\omega}.
$$
\end{lem}
\noindent {\bf Proof.} See [RZ, Proposition 9.1.11, p.
370].$\hfill\Box$

\medskip

Our next two lemmas deal with a special case for which the
principal result holds.
\begin{lem}
 Suppose that the set of positive integers $t$ for which there exists
 a (continuous) epimorphism
 $$
 C(\Gamma)\longrightarrow\hat{F}_t
 $$
 is not bounded. Then the principal result holds.
\end{lem}
\noindent {\bf Proof.} This follows from the proof of [Za1, Lemma
4.6].$\hfill\Box$\\

\smallskip

An immediate application is the following.
\begin{lem}
Suppose that, for all $\mathfrak{q}$, the closure of
$U(\mathfrak{q})$ in $\:\overline{\Gamma}$, $
\overline{U}(\mathfrak{q})$, is open in $\overline{\Gamma}$. Then
the principal result holds.
\end{lem}
\noindent{\bf Proof.}  The hypothesis ensures that
$|\overline{\Gamma}(\mathfrak{q}): \overline{U}(\mathfrak{q})|$ is
finite. We confine our attention to those (all but finitely many)
$\mathfrak{q}$ for which Theorem 3.5 and Lemma 5.1 hold. Let
$C(\Gamma)=C$.
 \noindent Now $C \cdot
\hat{U}(\mathfrak{q})$ is of finite index in $C \cdot
\hat{\Gamma}(\mathfrak{q}) = \hat{\Gamma}(\mathfrak{q})$. It
follows that
$$
C/C \cap \hat{U}(\mathfrak{q}) \cong C \cdot
\hat{U}(\mathfrak{q})/ \hat{U}(\mathfrak{q})
$$
is an open subgroup of
$$
\hat{\Gamma}(\mathfrak{q})/\hat{U}(\mathfrak{q})\cong\hat{F}_r.
$$
By [RZ, Corollary 3.6.4, p.~119] $C$ maps {\it onto}
$\hat{F}_{r'}$, for some $r'\geq r=r(\mathfrak{q})$. The result
follows from Lemmas 3.7 and 5.6.$\hfill\Box$

\medskip

Lemma 5.7 applies, for example, to the case
$\mathbf{G}=\mathbf{SL}_2, \;S=\{v\}$ and $\Gamma = {\rm
SL}_2(\mathcal{O})$ (as demonstrated in [Za1]). It is known [Mas1,
Theorem 3.1] that, when $\Gamma = {\rm SL}_2(\mathcal{O})$, the
``smallest congruence subgroup" of $\Gamma$ containing
$U(\mathfrak{q})$,
$$
\bigcap_{\mathfrak{q}' \neq \{0\}}
U(\mathfrak{q})\cdot\Gamma(\mathfrak{q}')=\Gamma(\mathfrak{q}),
$$
for all $\mathfrak{q}$.  It follows that in this case
$\overline{\Gamma}(\mathfrak{q}) = \overline{U}(\mathfrak{q})$,
for all $\mathfrak{q}$.

\medskip

 We now make use of the Strong Approximation Property for
$\mathbf{G}$. We will identify
$\overline{\mathbf{G}}(\mathcal{O})$ with the restricted
topological product $ \mathbf{G}(\hat{\mathcal{O}}) $. (See
Section $1$.) We record a well-known property.

\begin{lem}
For all $v\not\in S$, $\mathbf{G}(\mathcal{O}_v)$ is virtually a
pro-$p$ group.
\end{lem}
\noindent {\bf Proof.} In the notation of Section 1, the subgroup
$\mathbf{G}(\mathfrak{m})$ is of finite index in
$\mathbf{G}(\mathcal{O}_v)$ and is a pro-$p$ group. (See, for
example, [PR, Lemma 3.8, p.~138].)$\hfill\Box$

\medskip

It is convenient at this point to simplify our notation. We put
$$
C=C(\Gamma)\; \mbox{ and }\; \Lambda=\Gamma(\mathfrak{q}).
$$
It will always be assumed that Theorem 3.5 applies to
$\mathfrak{q}$ and (by Lemma 3.7) that $r(\mathfrak{q})\geq 2$. We
identify $\overline{\Lambda}$ with its embedding in
$\mathbf{G}(\hat{\mathcal{O}})$, (via the "diagonal" embedding of
$\Lambda$). We also identify each $\mathbf{G}(\mathcal{O}_v)$ with
its embedding as a normal subgroup of
$\mathbf{G}(\hat{\mathcal{O}_v})$. Let
$$
\phi: \hat{\Lambda}\longrightarrow \overline{\Lambda}
$$
denote the natural epimorphism.

\begin{lem}
For each $v\not\in S$, the group
$N_v:=\phi^{-1}(\overline{\Lambda}\cap\mathbf{G}(\mathcal{O}_v))$
is a closed, normal subgroup of $\hat{\Lambda}$ containing $C$.
Moreover, $$ N_v\cong P_v\amalg N_v(p), $$ where
\begin{itemize}
\item[(i)] $P_v$ is a projective group, all of whose proper, open
subgroups are isomorphic to $\hat{F}_{\omega}$;
\item[(ii)] $N_v(p)$ is a closed subgroup of $\hat{U}(\mathfrak{q})$ and is a free
profinite product of nilpotent
pro-$p$ groups, each of class at most $2$ and each generated by
torsion elements of $p$-power order.
\end{itemize}
\end{lem}
\noindent{\bf Proof.} Follows from Lemma 5.2.$\hfill\Box$

\medskip

Our next lemmas will be used to establish another condition under
which the principal result holds.
\begin{lem}
Let $\mid\mathbf{G}(\mathcal{O}):\Lambda\mid= n$ and let
$$
\pi(\overline{\Lambda})\,:=\,\prod_{v\not\in
S}(\overline{\Lambda}\cap\mathbf{G}(\mathcal{O}_v)).
$$
Then $ g^{n!} \in \pi(\overline{\Lambda})$ for all $g\in
\mathbf{G}(\hat{\mathcal{O}}).$
\end{lem}
\noindent{\bf Proof.} Since
$$
\mid\mathbf{G}(\mathcal{O}_v):\overline{\Lambda}\cap\mathbf{G}(\mathcal{O}_v)
\mid=\mid\overline{\Lambda}\cdot\mathbf{G}(\mathcal{O}_v):
\overline{\Lambda}\mid\leq\mid\mathbf{G}(\hat{\mathcal{O}}):\overline{\Lambda}\mid\leq
n,
$$
the assertion follows. $\hfill\Box$

\begin{lem}
With the above notation,
$$
\mid\mathbf{G}(\hat{\mathcal{O}}):\pi(\overline{\Lambda})\;.\;\overline{U}(\mathfrak{q})\mid<\infty.
$$
\end{lem}
\noindent{\bf Proof.} Set $
\Lambda^*:=\,\overline{\Lambda}/(\pi(\overline{\Lambda})\cdot\overline{U}(\mathfrak{q})).
$ The (compact, Hausdorff) group $\Lambda^*$ is finitely generated
by Lemma 5.1 and periodic by Lemma 5.10. It follows from
Zel'manov's celebrated result [Ze] that $\Lambda^*$ is finite.
$\hfill\Box$

\medskip

We are now able to prove the principal result.
\begin{thm}If $\Gamma$ is non-uniform, then
$$
C(\Gamma)\cong\hat{F}_{\omega}\amalg N(\Gamma),
$$
where $N(\Gamma)$ is a free profinite product of elementary
abelian pro-$p$ groups, each isomorphic to the direct product of
$2^{{\aleph}_0}$ copies of $\mathbb{Z}/p\mathbb{Z}$.
\end{thm}
\noindent{\bf Proof.} There are two possibilities, the first of
which can be readily dealt with.

\medskip

\noindent{\bf Case A:} {\it For all} $\mathfrak{q}$,  {\it we
have} $P_v\leq C$ {\it , for all} $v\notin S$.

 \medskip

\noindent
 For all $\mathfrak{q}$ and all $v\not\in S$, it follows from Lemma
 5.9
 that
 $
 \pi(\overline{\Lambda})\leq\overline{U}(\mathfrak{q}).
 $
 The principal result then follows from Lemmas 5.7 and 5.11. We
 consider the remaining case.

 \medskip

 \noindent {\bf Case B:} {\it There exists} $\mathfrak{q}$ {\it
 and}
 $v\not\in S$ {\it such that} $P_v\nleq C$.

 \medskip

\noindent
 For such a $v$ there exists an open, normal subgroup $L$ of $N_v$,
 containing $C$, such that
 $
 L\cap P_v\neq P_v.
 $
 It follows from Lemma 5.4 that
 $$
 L\cong \hat{F}_{\omega}\amalg\cdots\,.
 $$
 Restricting $\phi$ to $L$, there are again two possibilities. If
 $\phi(\hat{F}_{\omega})$ is trivial, then
 $C\cap\hat{F}_{\omega}=\hat{F}_{\omega}$. Since $C$ is a closed
 normal subgroup of $L$, the principal result follows from Lemmas
 5.4 and 5.5.\\

 Thus, from now we may assume that
 $\phi(\hat{F}_{\omega})$ is non-trivial. Note that
 $$
 L\cong \hat{F}_n\amalg \cdots\,
 $$ for all $n\ge 2$.
 Again restricting $\phi$, to $L$ there are two cases.

 \medskip

 \noindent {\bf Subcase B (i)}: $\phi(\hat{F}_n)$ {\it is finite for all} $n\geq
 2$.

 \medskip

\noindent
 It follows that, for all $n\geq 2$ we have that
 $
 C\cap\hat{F}_n\cong\hat{F}_{n'}
 $
 for some $n'\geq n$; see [RZ, Theorem 3.6.2, p.~118]. Then, as $C$ is a closed, normal subgroup of
 $L$,
 $$
 C\cong\hat{F}_{n'}\amalg\cdots\,
$$
by Lemma 5.4. Thus $C$ maps onto $\hat{F}_{n'}$. The principal
result follows from Lemma 5.6.

\medskip

\noindent {\bf Subcase B (ii)}: {\it There exists} $n\geq 2$ {\it
such that} $\phi(\hat{F}_n)$ {\it is infinite}.

\medskip

\noindent We consider $\phi(\hat{F}_n)$ as a subgroup of
$\mathbf{G}(\mathcal{O}_v)$. Let $M=\mathbf{G}(\mathfrak{m})$, as
defined in the proof of Lemma 5.8. Then
$$
(\phi^{-1}(M\cap\phi(\hat{F}_n)))\cap\hat{F}_n\cong\hat{F}_{n'}
$$
for some $n'\geq n$, by [RZ, Theorem 3.6.2, p.~118], and,
intersecting both sides with $C$, it follows that
$$
C\cap\hat{F}_n=\,C\cap \hat{F}_{n'}.
$$
Suppose that $M\cap\phi(\hat{F}_n)$ is non-abelian. Then by [BL]
and Lemma 5.8 this group is not free pro-$p$ and hence does not
satisfy Schreier's formula [RZ, p.~320], by [RZ, Theorem 8.4.7,
p.~321]. It follows that $\hat{F}_n/C\cap\hat{F}_n$ does not
satisfy Schreier's formula. But then
$$
C\cap\hat{F}_n\cong\hat{F}_{\omega}
$$
thanks to [RZ, Corollary 8.4.4, p.~320]. The principal result
follows from Lemmas 5.4 and 5.5.

\medskip

It remains to consider the case where $M\cap\phi(\hat{F}_n)$ is a
finitely generated, infinite abelian group. Then by [RZ, Lemma
8.4.5, p.~320] this group does not satisfy the Schreier formula
(in which case the principal result holds as above)  {\it unless}
it is infinite cyclic. In the latter case we can use [RZ, Theorem
8.4.3, p.~319] to conclude that again
$$
C\cap\hat{F}_n\cong\hat{F}_{\omega},
$$
from which the principal result follows as above. $\hfill\Box$\\

\noindent {\bf Acknowledgement.} The authors are grateful to
Professor Gopal Prasad for clarifying a number of points, in
particular for providing the reference in [BT2] used in Lemma 3.4.\\
\\
\noindent {\bf Added comment.} The referee has indicated that a
recent paper of P. Gille, ``Unipotent subgroups of reductive groups
in characteristic $p>0$", Duke Math. J. 114 (2002), 307-328, can be
used to provide a shorter proof of Theorem C. However our proof is
more elementary. Moreover our proof can, in principle, be
generalized to the case of a reductive group of $K$-rank 1, where
$K$ is any field of positive characteristic. (See the Remark in
Section 4.) Gille's argument only applies to non-zero characteristic
fields with some extra special properties.

\end{document}